\documentclass[10pt,reqno, a4paper]{amsart}

\usepackage{amsmath,amsthm}
\usepackage{amssymb}
\usepackage{amsfonts}
\usepackage{amscd}
\usepackage{enumerate}
\usepackage{color}
 at 9truept

\newtheorem{thm}{Theorem}[section]
\newtheorem{theorem}[thm]{Theorem}

\newtheorem{lemma}[thm]{Lemma}
\newtheorem{proposition}[thm]{Proposition}
\newtheorem{definition}[thm]{Definition}

\theoremstyle{definition}
\newtheorem{example}[thm]{Example}
\newtheorem{examples}[thm]{Examples}
\newtheorem{remark}[thm]{Remark}
\newtheorem{remarks}[thm]{Remarks}

\newtheorem{free text}[thm]{}

\numberwithin{equation}{section}

\newcommand{\N} {\mathbb{N}}
\newcommand{\Z} {\mathbb{Z}}

\newcommand{\C} {\mathbb{C}}

\newcommand{\Der}{\mathbb{D}er}

\newcommand{{\Uop}}{{U^{\rm op}}}
\newcommand{{\Ae}}{{A^{\rm e}}}
\newcommand{{\Aop}}{{A^{\rm op}}}
\newcommand{{\Bop}}{{B^{\rm op}}}

\newcommand{{\op}}{{{\rm op}}}
\newcommand{{\coop}}{{{\rm coop}}}
\newcommand{{\sop}}{{*^{\rm op}}}

\newcommand{\lact}{{\,\raise1pt\hbox{$\scriptscriptstyle{\rhd}$}\,}}                  %
\newcommand{\ract}{{\,\raise1pt\hbox{$\scriptscriptstyle{\lhd}$}\,}}                  
\newcommand{\blact}{{\,\raise1pt\hbox{$\scriptscriptstyle{\blacktriangleright}$}\,}}  %
\newcommand{\bract}{{\,\raise1pt\hbox{$\scriptscriptstyle{\blacktriangleleft}$}\,}}   %

\begin{document}

\title{Differential calculus over double Lie algebroids}
\author{Sophie Chemla}

\maketitle

\begin{abstract}
M. Van den Bergh (\cite{VdB2}) defined the notion of a double Lie algebroid and showed  that a  double quasi-Poisson algebra gives rise to a double Lie algebroid. We give new examples of double Lie algebroids and develop a differential calculus in that context. We recover the non commutative Karoubi--de Rham complex (\cite{Karoubi},\cite{CB-PE-VG}) and the double Poisson--Lichnerowicz cohomology (\cite {Pichereau-vanWeyer}) as particular cases of our construction.
\end{abstract}
\maketitle

\section{Introduction}

Let $k$ be a field of characteristic 0. The pair 
 $({\mathcal L}, \omega)$ is a Lie--Rinehart algebra (\cite{Rinehart}) over a commutative $k$-algebra $A$ if ${\mathcal L}$ is endowed with a $k$-Lie algebra structure and an $A$-module structure where  the two structures are linked by a compatibility relation involving the anchor $\omega$. There is a 
 one-to-one correspondence (\cite{Kosmann}) between Lie--Rinehart  algebra structures on ${\mathcal L}$ and Gerstenhaber algebra structures on 
 $\Lambda_A {\mathcal L}$ (where $\Lambda_A {\mathcal L}$ is the exterior algebra of 
 ${\mathcal L}$).
 If $A$ is smooth and ${\mathcal L}$ is a finitely generated  projective $A$-module, then ${\mathcal L}$ is the $A$-module  of global sections of  a Lie algebroid. 
 Lie--Rinehart algebras generalize at the same time $A$-Lie algebras (case where $\omega$ is 0) and the Lie algebra of derivations over $A$ (then $\omega=id$). 
 A Poisson smooth algebra gives rise to a Lie--Rinehart algebra structure (on ${\mathcal L}=\Omega^1$)
 and this is also true for quasi-Poisson smooth $G$-algebras 
 (\cite{A-KS-M}, \cite{B-C}, if $X$ is a Poisson manifold, then 
 ${\mathcal L}=\Omega^1_X \oplus \Gamma (X\times {\mathfrak g}))$. 
 Lie algebroids have been extensively studied and used, in particular in Poisson geometry. 
 A Lie--Rinehart algebra ${\mathcal L} $ defines a differential $d_{\mathcal L}$ on the graded algebra 
 $\Lambda _A {\mathcal L}^*$. In the case of an $A$-Lie algebra, we recover the Cartan--Eilenberg differential. In the case where  ${\mathcal L}=Der (A)$ (with $A$ smooth), we recover the de Rham differential. In the case where ${\mathcal L}=\Omega^1(X)$ ($X$ being a Poisson manifold), we recover the Lichnerowicz--Poisson differential. More generally a differential calculus has been  developed for Lie algebroids (\cite{E-L-W}).
 
 In this article, we are interested in the case where $A$ is not necessarily commutative. 
 As in \cite{VdB1} and \cite{VdB2}, we use a non commutative version of differential geometry based on an idea of Kontsevitch : For a property of a non commutative $k$-algebra to have a geometric meaning, it should induce (through a trace map) a standard geometric property on all representation spaces $Rep(A,N)=Hom(A, Mat_N(k))$, for all integer $N$. The coordinate ring of $Rep (A,N)$ is 
 $${\mathcal O}_{Rep(A,N)}=
 \frac{k[a_{i,j} \mid a \in A, \; (i,j)\in [1,N]^2]}{<(ab)_{i,j}=a_{i,k}b_{k,j}, \quad a,b \in A>}.$$
 If the  $k$-algebra $A$ is smooth, then  the commutative $k$-algebra
 ${\mathcal O}_{Rep(A,N)}$ is smooth for all $N$ (\cite{VdB1}, \cite{VdB2}). 
 If $A=k<x_1, \dots , x_m>$ is a free associative $k$-algebra generated by the variables $x_1, \dots , x_m$,  then 
 $$Rep_N(A)=Mat_N(k)\oplus Mat_N(k)\oplus \dots \oplus Mat_N(k)\simeq k^{mN^2}$$
where, if $\alpha \in [1,m], \quad  M_\alpha$ is  the $N\times N$ matrix $(x_{i,\alpha}^j)_{i,j}$ and 
$Tr (x_{\alpha_1}\dots x_{\alpha_r})=Tr(M_{\alpha_1}M_{\alpha_2}\dots M_{\alpha_r})$. 
 
 The non commutative geometry notions defined according to the  principle of Kontsevitch are often  called by a name  with the prefix "double".  Thus derivations  $Der (A)$ are replaced by double derivations 
 $\mathbb{D}er (A)$ (that are  derivations from $A$ to $A\otimes A$ considered with its exterior $A$-bimodule structure), Poisson algebras are replaced by  double Poisson algebras, etc ...
 The Karoubi--de Rham complex  was defined in the non commutative setting (\cite{Karoubi}), 
 the contraction  and the Lie derivative by a double vector field were defined in  (\cite{CB-PE-VG}).
 Double Poisson cohomology was  defined (\cite{vanWeyer}, \cite{Pichereau-vanWeyer}) and   computed  for some examples of double Poisson brackets associated with quivers.

 In \cite{VdB2},  double Lie algebroids or double Lie--Rinehart algebras were   introduced   and it was shown that a double quasi-Poisson algebra gives rise to a double Lie algebroid. A double Lie  algebroid $\mathbb L$ being an $A$-bimodule, one can construct the tensor algebra of 
 $\mathbb{L}$, 
 $T_A({\mathbb L})$. In \cite{VdB2}, 
  the double Lie  algebroid structure has been defined by the 
 Gerstenhaber double structure it induces on $T_A({\mathbb L})$. 
 
 We give a direct definition of a double Lie algebroid  and show that there is a correspondence between double Lie  algebroid structures on an $A$-bimodule ${\mathbb L}$ and double Gerstenhaber algebra structures on $T_A({\mathbb L})$. 
 We give new examples  of double Lie  algebroids. Studying the  case where $A=k$, we see that any associative $k$-algebra has a natural double Lie algebroid structure. Then, 
 we develop a differential calculus for double Lie  algebroids : definition of a differential $d_{\mathbb L}$, for which we give an explicit formula, Lie derivative, contraction, etc... In the case where 
 ${\mathbb L}={\mathbb D}er (A)$ and $A$ is smooth, we recover the Karoubi--de Rham complex 
 (\cite{Karoubi}, \cite{CB-PE-VG}) but some of our formulas are new even in that case. 
 In the case where 
 $A$ is a smooth double Poisson algebra and ${\mathbb L}=\Omega_A^1$, the differential 
 $d_{\mathbb{L}}$ coincides with that of  double Poisson cohomology 
 (\cite{vanWeyer},\cite{Pichereau-vanWeyer}).  
 In the case where $A=k$, $d_{\mathbb{L}}$ is the differential computing  cyclic cohomology. 
 Thus we recover the classical picture.   The theory  of double Lie algebroids encompasses 
several theories. \\

{\bf Acknowledgement}

I am most grateful to Y. Kosmann-Schwarzbach  
and V. Rubtsov for their help and encouragement, which were essential for the completion of this article. They both gave me several references that were instrumental for me to complete this article. They explained to me several aspects of the theory of double Poisson algebras and commented on earlier versions of this preprint. I am happy to extend to them my most sincere thanks. \\

{\bf Convention :} 

We will use the same notation as \cite{VdB1}

If $(V_i)_{i=1, \dots ,n}$ are $k$ vector spaces and $s \in S_n$, then for 
$a=a_1 \otimes \dots \otimes a_n \in V_1\otimes \dots \otimes V_n$
$$\tau_s(a)=a_{s^{-1}(1)}\otimes \dots \otimes a_{s^{-1}(n)}.$$
If $(V_i)_{i=1, \dots ,n}$ are $k$ graded vector spaces and $s \in S_n$, then for 
$a=a_1 \otimes \dots \otimes a_n \in V_1\otimes \dots \otimes V_n$
$$\sigma_s(a)=(-1)^ta_{s^{-1}(1)}\otimes \dots \otimes a_{s^{-1}(n)}.$$
where $t={\displaystyle \sum_{i<j, s^{-1}(i)>s^{-1}(j)}} \vert a_{s^{-1}(i)}\vert \vert a_{s^{-1}(j)}\vert$.

$\tau_{(12)}$ (respectively $\sigma_{(12)}$) will be also denoted $(x\otimes y)^\circ =y\otimes x$
 (respectively $(x\otimes y)^\circ=(-1)^{\vert x\vert \vert y\vert} y \otimes x$.)

Let $B$ be a fixed $k$-algebra that will be, most of the time, semisimple commutative of the form 
$B=ke_1\oplus \dots \oplus ke_n$ with $e_i^2=e_i$. A $B$-algebra is a $k$-algebra $A$ equipped with a morphism of $k$-algebras $B \to A$. The notion of $B$-algebra allows to define  relative versions.

\section{Definitions and generalities}

Most of the definition and results of this section come from  \cite{VdB1}.

\begin{definition} An $n$-bracket is a linear map 
$$\{\!\{-, \dots ,-\}\!\} : A^{\otimes n} \to A^{\otimes n}$$
which is a derivation $A \to A^{\otimes n}$ in its last argument for the outer bimodule structure on 
$A^{\otimes n}$ i. e
$$\{\!\{a_1, a_2, \dots , a_{n-1},a_n a_n^\prime \}\!\}=
a_n \{\!\{a_1, a_2, \dots , a_{n-1}, a_n^\prime \}\!\}+ \{\!\{a_1, a_2, \dots , a_{n-1},a_n \}\!\} a_n^\prime $$
and which is cyclically antisymmetric in the sense that 
$$\tau_{(1\dots n)} \circ \{\!\{-,\dots, -\}\!\} \circ \tau_{(1\dots n)}^{-1}=
(-1)^{n+1}\{\!\{-, \dots , -\}\!\}.$$
If $A$ is a $B$-algebra, then an $n$-bracket is $B$-linear if it vanishes when its last argument is in the image of $B$. 
\end{definition}

As in \cite{VdB1},  we set 
$$\{\!\{a,b\}\!\}_L=\{\!\{a,b_1\}\!\} \otimes b_2\otimes \dots \otimes b_n$$
Associated to a double bracket $\{\!\{-,-\}\!\}$, we define a tri-ary operation 
$\{\!\{-,-,-\}\!\}$ as follows :
$$\{\!\{a,b,c\}\!\}= \{\!\{a, \{\!\{b,c\}\!\} \}\!\}_L + 
\tau_{(123)}\{\!\{b, \{\!\{c,a\}\!\} \}\!\}_L + \tau_{(132)}\{\!\{c, \{\!\{a,b\}\!\} \}\!\}_L $$
\begin{proposition}(\cite{VdB1})
$\{\!\{-,-,-\}\!\}$ is a 3-bracket.
\end{proposition}
\begin{definition}(\cite{VdB1}). 
Let $A$ be a $k$-algebra. A double bracket $\{\!\{-,-\}\!\}$ on A is a double Poisson bracket if 
$\{\!\{-,-,-\}\!\}=0$. An algebra with a double Poisson bracket is a double Poisson algebra. 
\end{definition}

\begin{example}  (\cite{VdB1}, \cite{Powell}) One may characterize the double Poisson brackets on $k[t]$.  For $\lambda , \mu , \nu \in k$,
$$\{\!\{t,t\}\!\}=\lambda (t\otimes 1 -1 \otimes t) +\mu (t^2\otimes 1 -1 \otimes t^2) +\nu (t^2 \otimes t-t \otimes t^2)$$
defines a double Poisson structure if and only if $\lambda \nu -\mu^2=0$ and any double Poisson structure on $k[t]$ is of this form.  

We will see many other examples of double Poisson algebras further. 

\end{example}

The following proposition was proved in \cite{VdB1} :
\begin{proposition}
Assume that $(A, \{\!\{-,-\}\!\})$ is a double Poisson algebra.  For any elements $a$ and $b$ of $A$, we set $\{a,b\}=\{\!\{a,b\}\!\}^\prime \{\!\{a,b\}\!\}^{\prime \prime}.$
Then the following holds ~:

(1) $\{-,-\}$ is a derivation in its second argument and vanishes on commutators in its first argument.

(2) $\{-,-\}$ is anti-symmetric modulo commutators.

(3) $\{-,-\}$ makes $A$ into a left Loday algebra,  i.e $\{-,-\}$ satisfies the following version of Jacobi identity 
$$\{a, \{b,c\}\}=\{\{a,b\},c\}+ \{b, \{a,c\}\}$$

(4) $\{-,-\}$ makes $A/[A,A]$ into a Lie algebra.
\end{proposition}

\begin{definition}(\cite{CB}, \cite{CB1})
A  Poisson structure on $A$ is a Lie bracket $\{-,-\}$ on $\dfrac{A}{[A,A]}$ such that for each $a\in A$, the map 
$\{ \overline{a},-\} : \dfrac{A}{[A,A]} \to \dfrac{A}{[A,A]}$ is induced by a derivation on $A$. 
\end{definition}
\begin{remarks}\label{$H_0$-Poisson}
\begin{enumerate}

\item In the case where $A$ is commutative, we recover the usual Poisson bracket.

\item It was shown in \cite{VdB1} (lemme 2.6.2) that a double Poisson bracket on $A$ induces a Poisson structure  on $A$. 

\item A Poisson structure on $A$ is also called a $H_0$-Poisson structure as it is a structure on 
$\dfrac{A}{[A,A]}=HC_0(A)$ which is the $0$th cyclic group of $A$. In \cite{B-C-E-R}, derived Poisson structure were defined on higher cyclic cohomology group. \\
\end{enumerate}
\end{remarks}
Let $D$ be a graded algebra. There are two commuting $D^e$-module structures on $D \otimes D$ : For any homogeneous elements $\alpha, \beta, x, y$ in $D$, 
$$\begin{array}{l}
\alpha (x\otimes y )\beta=\alpha x\otimes y \beta\\
\alpha *(x\otimes y)*\beta=(-1)^{\vert \alpha \vert \vert \beta\vert + \vert \alpha \vert \vert x\vert +\vert y \vert \vert \beta \vert }
x\beta \otimes  \alpha y\\
\end{array}$$
\begin{definition}(\cite{M-T1})
Let $d \in \Z$ and let  $D$ be a graded algebra.  $D$ is called a double Gerstenhaber algebra of degree $d$ if it is equipped with a graded bilinear map 
$$\{\!\{-,-\}\!\} : D \otimes D \to D \otimes D$$
of degree d such that the following identities hold :
$$\begin{array}{l}
1) \{\!\{\alpha,\beta \gamma \}\!\}=
(-1)^{(\vert \alpha \vert +d)\vert \beta \vert }\beta \{\!\{ \alpha,\gamma \}\!\} + \{\!\{\alpha,\beta\}\!\} \gamma\\
1^\prime ) \{\!\{\beta \gamma , \alpha\}\!\}=
(-1)^{(\vert \alpha \vert +d)\vert \beta \vert }\beta* \{\!\{ \alpha,\gamma \}\!\} +
 \{\!\{\alpha ,\beta \}\!\}* \gamma \\
2) \{\!\{\alpha, \beta\}\!\}=-(-1)^{(\vert \alpha \vert +d)(\vert \beta \vert +d)}
\sigma_{(12)}\{\!\{\beta, \alpha \}\!\}\\
3) \{\!\{\alpha, \{\!\{\beta,\gamma\}\!\} \}\!\}_L + 
(-1)^{(\vert \alpha \vert +d)(\vert \beta \vert +\vert \gamma\vert )}
\sigma_{(123)}\{\!\{\beta, \{\!\{\gamma,\alpha\}\!\} \}\!\}_L
+(-1)^{(\vert \gamma \vert +d)(\vert \alpha \vert + \vert \beta\vert )}
\sigma_{(132)}\{\!\{\gamma , \{\!\{\alpha ,\beta \}\!\} \}\!\}_L=0
\end{array}$$
\end{definition}
\begin{remarks} 1)  The definition of double Gerstenhaber algebra is given in \cite{VdB1}. It is extended to the case of 
double Gerstenhaber algebra of degree $d\in \Z$ in \cite{M-T1}. The case $d=-1$ corresponds to double  Gerstenhaber algebras (see \cite{VdB1}) and the case $d=0$ corresponds to double Poisson algebras (\cite{VdB1}).

2) Assertions $1)$ and $1^\prime)$ are equivalent if assertion 2) is satisfied (\cite{M-T1}).\\
\end{remarks}
If $D$ is a double Gerstenhaber algebra, we define the associated bracket $\{-,- \} : D \otimes D \to D$ by :
$$\forall (\alpha, \beta)\in D^2, \quad  \{\alpha, \beta\}=
\{\!\{\alpha, \beta\}\!\}^\prime \{\!\{\alpha ,\beta\}\!\}^{\prime \prime}.$$
The following proposition was partly stated in \cite{VdB1} : 

\begin{proposition}\label{properties of d-double Gerstenhaber algebras}
Let $D$ be a Gerstenhaber algebra of degree $d$ and let $\alpha , \beta , \gamma$ three homogeneous elements in $D$.

1) $\{ \alpha \otimes \beta - (-1)^{\vert \alpha \vert \vert \beta \vert}\beta \otimes \alpha, \gamma \}=0.$

2) $\{ \alpha , \beta \}- (-1)^{(\mid \alpha \mid +d)(\mid \beta \mid +d)} \{\beta , \alpha \} \in [D, D]$

3) $\{ \alpha , \{\!\{ \beta , \gamma \}\!\}\}-\{\!\{ \{ \alpha ,\beta \}, \gamma \}\!\} -
(-1)^{(\mid \alpha \mid +d)(\vert \beta \vert +d)}
\{\!\{\beta, \{\alpha , \gamma\} \}\!\}=0$

where $\{\alpha, -\}$ acts on tensors by 
$$\{ \alpha , \{\!\{ \beta , \gamma \}\!\} \}=\{ \alpha , \{\!\{ \beta , \gamma \}\!\}^\prime \} \otimes \{\!\{ \beta , \gamma \}\!\}^{\prime \prime}+
(-1)^{(\vert \alpha \vert +d) \vert \{\!\{ \beta, \gamma\}\!\}^\prime \vert} 
\{\!\{ \beta , \gamma \}\!\}^\prime  \otimes  \{\alpha , \{\!\{ \beta , \gamma \}\!\}^{\prime \prime } \} $$

4) $\dfrac{D}{[D,D]}[d]$ is a  graded Lie algebra. 
\end{proposition}

\begin{remark}
Proposition \ref{properties of d-double Gerstenhaber algebras} is proved in the ungraded case in \cite{VdB1}. In this case,  $d=0$ and $D$ is a   double Poisson algebra. Our proof is similar to that of \cite{VdB1} so that we only sketch the main lines of it .
\end{remark}

\begin{proof}
1) and 2)  are  straightforward computations if one uses the relation 
$$\{\!\{\alpha, \beta \}\!\}=
-(-1)^{(\vert \alpha \vert  +d)(\vert \beta \vert +d)+
\vert \{\!\{\alpha, \beta\}\!\}^\prime\vert \vert \{\!\{\alpha,\beta\}\!\}^{\prime \prime}\vert}
\tau_{(12)}\{\!\{\beta,\alpha\}\!\}$$

4) is a consequence of  the previous statements. Let us now prove 3). 

We will make use of the following lemma whose proof is left to the reader:

\begin{lemma} Set 
$\{\!\{ \alpha , \{\!\{ \gamma , \beta \}\!\} \}\!\}_R=
(-1)^{\mid \alpha \mid +d)\mid \{\!\{ \gamma ,\beta \}\!\}^\prime\vert }
 \{\!\{ \gamma , \beta \}\!\}^\prime \otimes  
\{\!\{ \alpha ,\{\!\{ \gamma , \beta \}\!\}^{\prime \prime} \}\!\}$. 
The following equality holds 
$$\{\!\{ \alpha , \{\!\{ \gamma , \beta \}\!\} \}\!\}_R=- (-1)^{(\vert \beta \vert +d )(\vert \gamma \vert +d)}
{\color{black}\sigma}_{(123)} \{\!\{ \alpha , \{\!\{ \beta , \gamma\}\!\} \}\!\}_L\\ 
$$
\end{lemma}

Let us now compute the three terms of the equality 3). Assertion 3) will follow from these computations. 

$$\begin{array}{rcl}
\{\alpha , \{\!\{ \beta , \gamma \}\!\} \}&=&
\{ \alpha , \{\!\{ \beta , \gamma \}\!\}^{\prime}\} \{\!\{ \beta , \gamma \}\!\}^{\prime \prime}
+ (-1)^{(\vert \alpha \vert +d) \vert \{\!\{ \beta , \gamma \}\!\}^\prime \vert } 
\{\!\{ \beta , \gamma \}\!\}^{\prime}
\{ \alpha , \{\!\{ \ \beta , \gamma \}\!\}^{\prime \prime}\}\\
&=& (m\otimes 1) \{\!\{ \alpha, \{\!\{\beta , \gamma \}\!\} \}\!\}_L + 
(1 \otimes m) \{\!\{ \alpha, \{\!\{ \beta , \gamma \}\!\} \}\!\}_R\\
&=& (m\otimes 1) \{\!\{ \alpha, \{\!\{\beta , \gamma \}\!\} \}\!\}_L - 
(1 \otimes m) {\color{black} \sigma}_{(123)}\{\!\{ \alpha, \{\!\{ \gamma , \beta \}\!\} \}\!\}_L 
(-1)^{(\vert \beta \vert +d)(\vert \gamma \vert +d)}\\
\end{array}$$

$$\begin{array}{l}
\{\!\{ \{ \alpha, \beta \}, \gamma \}\!\}\\
=-{\color{black} \sigma}_{(12)}\{\!\{ \gamma, \{\alpha ,\beta\} \}\!\}
(-1)^{(\vert \gamma \vert +d) (\vert \alpha \vert +\vert \beta \vert )}\\
=-{\color{black}\sigma}_{(12)}\left [ \{\!\{ \alpha , \beta \}\!\}^\prime \{\!\{\gamma, \{\!\{\alpha , \beta \}\!\}^{\prime \prime}\}\!\}\right ]
(-1)^{(\vert \gamma \vert +d)(\vert \alpha \vert +\vert \beta \vert )} -
{\color{black}\sigma}_{(12)}\left [ \{\!\{\gamma , \{\!\{ \alpha , \beta \}\!\}^{\prime  }\}\!\} \{\!\{ \alpha , \beta \}\!\}^{\prime \prime }\right ]
(-1)^{(\vert \gamma \vert +d)(\vert \alpha \vert +\vert \beta \vert )}\\
= -(1 \otimes m) {\color{black}\sigma}_{(132)}\{\!\{\gamma , \{\!\{ \alpha , \beta \}\!\} \}\!\}_R
(-1)^{(\vert \gamma \vert +d)(\vert \alpha \vert +\vert \beta \vert )}
-(m\otimes 1) {\color{black}\sigma_{(132)}}\{\!\{\gamma , \{\!\{ \alpha , \beta \}\!\}\}\!\}_L
(-1)^{(\vert \gamma \vert +d)(\vert \alpha \vert +\vert \beta \vert )}\\
= (1 \otimes m) {\color{black}\sigma}_{(132)}\{\!\{\gamma , \{\!\{ \beta , \alpha \}\!\}\}\!\}_L
(-1)^{(\vert \gamma \vert +d)(\vert \alpha \vert +\vert \beta \vert )}
(-1)^{(\vert \alpha \vert +d)(\vert \beta \vert +d )}\\
\quad -(m\otimes 1) {\color{black}\sigma}_{(132)}\{\!\{\gamma , \{\!\{ \alpha , \beta \}\!\}\}\!\}_L
(-1)^{(\vert \gamma \vert +d)(\vert \alpha \vert +\vert \beta \vert )}
\end{array}$$

$$\begin{array}{rcl}
\{\!\{ \beta , \{ \alpha , \gamma\} \}\!\}&= &
\{\!\{ \beta , \{\!\{ \alpha , \gamma\}\!\}^\prime \}\!\}  \{\!\{ \alpha , \gamma\}\!\}^{\prime \prime}
+ (-1)^{(\vert \beta \vert +d) \mid \{\!\{ \alpha ,\gamma \}\!\}^\prime\mid }
\{\!\{ \alpha , \gamma\}\!\}^\prime  \{\!\{ \beta , \{\!\{ \alpha , \gamma\}\!\}^{\prime \prime} \}\!\}\\
&=& (1 \otimes m ) \{\!\{ \beta , \{\!\{ \alpha , \gamma\}\!\} \}\!\}_L 
+ (m\otimes 1) \{\!\{ \beta , \{\!\{ \alpha , \gamma\}\!\} \}\!\}_R\\
&=& (1 \otimes m ) \{\!\{ \beta , \{\!\{ \alpha , \gamma\}\!\} \}\!\}_L 
- (m\otimes 1) {\color{black}\sigma}_{(123)}\{\!\{ \beta , \{\!\{ \gamma , \alpha\}\!\} \}\!\}_L
(-1)^{(\vert \alpha \vert +d) (\vert \gamma \vert +d)}\\
\end{array}$$
\end{proof}

\section{Double derivations} 

In this section, we recall results of \cite{VdB1}. 

Let $B$ be a $k$-algebra. A $B$-algebra is a pair $(A, \eta )$ where $\eta : B \to A$ is an algebra morphism.

Denote by $m : A \otimes_B A \to A$ multiplication on $A$. One sets $\Omega^1_B(A):= \rm{Ker} (m)$.  It is naturally endowed with an $A^e$-module structure.  If $a \in A$, then 
$da=a\otimes 1 -1 \otimes a$ belongs to $\Omega^1_B(A)$. If $A$ is finitely generated as a $B$-algebra, then 
$\Omega^1_B A$ is finitely generated as a left $A^e$-module. \\

\begin{definition} A $B$-algebra $A$ is called smooth over $B$ if it is finitely generated as an $B$-algebra and 
$\Omega_B^1A$ is projective as a left $A^e$-module. \end{definition}

\begin{example}
Let ${\mathcal Q}=(Q, I)$ be a finite quiver with vertex set $I=\{1, \dots, n\}$ and edge set $Q$. Denote by $e_i$ the idempotent associated to the vertex $i$ and we put $B={\displaystyle \oplus_{i}}ke_i$. The path algebra  
$A=k[Q]$ is smooth over $B$ (\cite{VdB1}). 
\end{example}

Let  $Der_B (A, A \otimes A)$ be the $B$-derivations from $A$ to $A \otimes A$ where we put the outer bimodule structure on $A \otimes A$. Such a derivation is called a double derivation. 
Inner bimodule structure on $A^e$ allows to endow 
$Der_B (A, A \otimes A)$ with an $A^e$-module structure :
$$\alpha \cdot D \cdot \beta (a)=D(a)^\prime \beta \otimes \alpha D(a)^{\prime \prime}$$

{\bf Notation :} The $A^e$-module $Der_B (A, A \otimes A)$ will be denoted 
$\mathbb{D}er_B(A)$. 

\begin{example} We assume that we are in the situation where 
$B=ke_1 \oplus \dots \oplus ke_n$ with $e_i^2=e_i$. 

Let $E_i : A \to A \otimes A$ defined by 
$E_i(a)=ae_i \otimes e_i -e_i \otimes  {\color{black}e_ia}$ and 
$E=\displaystyle \sum_{i=1}^n E_i \in \mathbb{D}er_B (A)$.
\end{example}

\begin{proposition} (\cite{VdB1})
Let $\delta$, $\Delta$ in $\mathbb{D}er_B (A)$, then
$$\begin{array}{l}
\{\!\{\delta , \Delta \}\!\}_l^\sim = (\delta \otimes 1) \Delta - (1 \otimes \Delta)\delta \\
\{\!\{\delta , \Delta \}\!\}^\sim_r:= (1 \otimes \delta) \Delta - ( \Delta \otimes 1)\delta =
-\{\!\{\Delta , \delta \}\!\}^\sim_l \\
\end{array}$$
are derivation from $A$ to $A^{\otimes 3}$ where the $A^{e}$-bimodule structure on $A^{\otimes 3}$ is the outer structure. 
\end{proposition}
We define 
$$\begin{array}{rcl}
\{\!\{\delta ,\Delta\}\!\}_l=\tau_{(23)} \circ \{\!\{\delta , \Delta \}\!\}_l^\sim \in \mathbb{D}er (A) \otimes A\\
\{\!\{\delta , \Delta \}\!\}_r=\tau_{(12)} \circ \{\!\{\delta, \Delta\}\!\}_r^\sim \in A \otimes \mathbb{D}er(A)
\end{array}$$
We write 
$$\begin{array}{l}
\{\!\{\delta , \Delta \}\!\}_l= \{\!\{\delta , \Delta \}\!\}_l^\prime \otimes \{\!\{\delta , \Delta \}\!\}_l^{\prime \prime}\\
\{\!\{\delta , \Delta \}\!\}_r= \{\!\{\delta , \Delta \}\!\}_r^\prime \otimes \{\!\{\delta , \Delta \}\!\}_r^{\prime \prime}\\
\end{array}$$
with $\{\!\{\delta , \Delta \}\!\}_r^\prime , \{\!\{\delta , \Delta \}\!\}_l^{\prime \prime} \in A$ and 
$\{\!\{\delta , \Delta \}\!\}_l^\prime , \{\!\{\delta , \Delta \}\!\}_r^{\prime \prime} \in \mathbb{D}er (A)$. 
\begin{theorem} (\cite{VdB1})
For $a,b \in A$ and $\delta , \Delta \in \mathbb{D}er_B(A)$, the following equations 
$$\begin{array}{l}
\{\!\{a,b\}\!\}=0\\
\{\!\{\delta , a \}\!\}=\delta (a) \in A \otimes A\\
\{\!\{\delta , \Delta \}\!\}=\{\!\{\delta , \Delta\}\!\}_l + \{\!\{ \delta , \Delta \}\!\}_r
\end{array}$$
define a unique structure of double Gerstenhaber algebra on $T_A\mathbb{D}er_B(A)$.
\end{theorem}

{\bf Notation :} From now on, $T_A\mathbb{D}er_B(A)$ will be denoted $D_B(A)$. \\

\begin{proposition} \label{n-brackets}(\cite{VdB1})
Assume that $A$ is a finitely generated $k$-algebra. The linear map
$$\begin{array}{rcl}
\mu : \quad \quad (D_BA)_n & \to & \{$B$-linear \;$n$-brackets \;on \;$A$\}\\
Q=\delta_1 \dots \delta_n & \mapsto & \{\!\{-, \dots, - \}\!\}_Q=
\{\!\{-,\dots,-\}\!\}= {\displaystyle \sum_{i=0}^{n-1} (-1)^{(n-1)i}} \tau_{(1\dots n)}^i \circ \{\!\{-, \dots,-\}\!\}^\sim_Q \circ \tau_{(1\dots n)}^{-i}
\end{array}$$
where 
$$\{\!\{a_1, \dots ,a_n\}\!\}^\sim _{\delta_1\dots \delta_n} = 
\delta_n(a_n)^\prime \delta_1(a_1)^{\prime \prime} \otimes 
\delta_1(a_1)^\prime \delta_2(a_2)^{\prime \prime} \otimes \dots \otimes 
\delta_{n-1}(a_{n-1})^\prime \delta_n(a_n)^{\prime \prime}$$
This maps factors through $\dfrac{D_BA}{[D_BA,D_BA]}$. The map $\mu$ is an isomorphism if $A$ is 
$B$-smooth.
\end{proposition}

In  \cite{VdB1}, the following expression is proved for $\{\!\{-, \dots ,- \}\!\}_Q$. 
\begin{proposition}
For $Q \in (D_BA)_n$, the following identity holds :
$$\{\!\{a_1, \dots ,a_n \}\!\}_Q=(-1)^{\frac{n(n-1)}{2}}
\{\!\{a_1, \dots, \{\!\{a_{n-1}, \{Q,a_n\}\}\!\}_L \dots \}\!\}_L$$
\end{proposition}

The isomorphism $\mu$ allows to characterize  double Poisson algebras on a smooth algebra :

\begin{proposition}
Let $A$ be a smooth $B$-algebra. Double $B$-linear Poisson structures on $A$ are in bijection with the 
$P\in T^2\mathbb{D}er_B (A)$ such that $\{P,P\}=0$ modulo $[D_BA,D_BA]$.
\end{proposition}

\begin{remark}
In the case of a double bivector field, the formula of proposition \ref{n-brackets} is the double version of the formula $\{f,g\}=-[f,[P,g]]$ (\cite{Kosmann1}) for the Poisson bracket. 
\end{remark}

The double Poisson - Lichnerowicz cohomology was defined in 
\cite{vanWeyer}, \cite{Pichereau-vanWeyer} :

\begin{definition} Let $A$ be a $k$-double Poisson algebra with Poisson double bracket defined by the Poisson bivector 
$P\in T^2\mathbb{D}er_k (A)$ such that $\{P,P\}=0$ modulo $[D_kA,D_kA]$. The 
double Poisson - Lichnerowicz cohomology of $A$ is the cohomology of the complex
 $\left ( \dfrac{D_kA}{[D_k(A), D_k(A)]}, \{P,-\}\right )$. 
\end{definition}

The concept of quasi-Poisson G-manifolds was introduced  in \cite{A-KS-M} and   Massuyeau and Turaev (\cite{M-T2}) gave an algebraic formulation of it.  Double quasi-Poisson algebras were defined in 
\cite{VdB1}. They also give rise to  $H_0$-Poisson structures.

\begin{definition}
We assume that we are in the situation where 
$B=ke_1 \oplus \dots \oplus ke_n$ with $e_i^2=e_i$. 

Let $E_i : A \to A \otimes A$ defined by $
E_i(a)=ae_i \otimes e_i -e_i \otimes {\color{black}e_ia}$ and 
$E=\displaystyle \sum_{i=1}^n E_i$.\\
A double quasi-Poisson brachet on $A$ is a $B$-linear bracket $\{\!\{,\}\!\}$ such that 
$$\{\!\{-,-,-\}\!\} = \{\!\{-,-,-\}\!\}_{E^3}.$$
\end{definition}

\begin{proposition} (\cite{VdB1}) 
Let $(A, \{\!\{-,-\}\!\}$ be a quasi-Poisson algebra. Then 

1) $(A,\{-,-\})$ is a left Loday algebra.

2) $(A,\{-,-\})$ induces a $H_0$-Poisson structure on A. 
\end{proposition}

\begin{definition} $P \in (D_BA)_2$ is a differential double quasi Poisson bracket if 
$$\{P,P\}=\dfrac{1}{6}E^3 \; mod \; [D_BA,D_BA]$$
\end{definition}

\begin{remark} If $A$ is smooth, then quasi Poisson bracket and differential quasi Poisson bracket are equivalent notions. \end{remark}

\begin{examples} 

1) Examples of Poisson double brackets and of quasi-Poisson double brackets over the path algebra of a the double of a quiver are given in \cite{VdB1}. 

2) Let $A$ be the free associative algebra on $n$ variables, $A=k<x_1, \dots ,x_n>$. Linear double Poisson structures are studied in \cite{Pichereau-vanWeyer}. Define 
$\{\!\{x_i,x_j\}\!\}=b^k_{i,j}x_k \otimes 1 -b_{jk}^k 1 \otimes x_k$ and extend this skew symmetric bracket to  a biderivation. This linear skew symmetric bracket is a double Poisson bracket if and only if $x_ix_j=b^k _{i,j}x_k$ is an associative multiplication.

3) Quadratic double Poisson brackets on $\C<X_1,X_2, \dots , X_n>$ have been studied  in 
\cite{O-R-S}.  
Define 
$\{\!\{x_i,x_j\}\!\}=r^{k,l}_{i,j}x_k \otimes x_l -
t^{{\color{black}k,l}}_{\color{black}i,j} x_kx_l \otimes 1 - 
t^{{\color{black}k,l}}_{\color{black}j,i} 1 \otimes x_kx_l$ and extend this bracket to a bi-derivation. For this quadratic double bracket to be Poisson, $r$ has to satisfy the associative Yang-Baxter equation. 

4) Given an oriented surface $\Sigma$ with boundary $\partial \Sigma$ and base point $*\in \partial \Sigma$, a quasi-Poisson double algebra structure is constructed on the group algebra 
$A=\mathbb{K}[\pi]$  of $\pi_1(\Sigma, *)$ in \cite{M-T2}. The Lie bracket on $\dfrac{A}{[A,A]}$ in this case is 
twice the Goldman Lie bracket.

5) Double Poisson structures on a semi-simple algebra (over an algebraically closed field) are described in  \cite{vanWeyer}.

\end{examples}

\section{Double Lie algebroids }

V. Roubtsov drew my attention to the fact that the 
 definition of a double Lie algebroid is given in \cite{VdB2}. But it is defined  by its characterization in terms of double Gerstenhaber algebras 
as explained in  proposition \ref{characterization Gerstenhaber} below. The definition we give is more explicit. 
\begin{definition}
A double Lie  algebroid is a quadruple $(\mathbb{L},A, \omega  , \{\!\{ -, -\}\!\})$ where 

\begin{itemize}
\item $\mathbb{L}$ is an $A^e$-module. 
\item 
$$\begin{array}{rcl}
\{\!\{,\}\!\} &:& \mathbb{L} \otimes \mathbb{L} \to   \mathbb{L} \otimes A \oplus A \otimes \mathbb{L} \\
&&(D, \Delta)  \mapsto \{\!\{D, \Delta \}\!\}^\prime \otimes \{\!\{D, \Delta \}\!\}^{\prime \prime}=
\{\!\{D, \Delta \}\!\}_l^\prime \otimes \{\!\{D, \Delta \}\!\}_l^{\prime \prime}+
\{\!\{D, \Delta \}\!\}_r^\prime \otimes \{\!\{D, \Delta \}\!\}_r^{\prime \prime} \\
&&( {\rm with}\quad  \{\!\{D, \Delta \}\!\}_l^\prime, \{\!\{D, \Delta \}\!\}_r^{\prime \prime}\in \mathbb{L}
\quad {\rm  and } \quad
\{\!\{D, \Delta \}\!\}_l^{\prime \prime}, \{\!\{D, \Delta \}\!\}_r^{ \prime}\in A)
\end{array}$$

 is a map satisfying 
  $\{\!\{ D, \Delta \} \!\}=-\{\!\{\Delta , D\}\!\}^{o}$.
\item $\omega :\mathbb{L} \to \Der_B (A)$ is a morphism of $A^e$-modules. 
\item $ \{\!\{\omega (D), \omega (\Delta )\}\!\}=\omega \left ( \{\!\{D, \Delta \}\!\} \right ) 
$ where we also denote by 
$\omega$ the extension of $\omega$ to a map from $A \otimes \mathbb{L} \oplus \mathbb{L} \otimes A$ to 
$A \otimes \mathbb{D}er_B (A)\oplus \mathbb{D}er_B(A) \otimes A$. 

\item Jacobi identity : If $D_1 , D_2 , D_3$ are elements of $\mathbb{L}$, one has :

$$\{\!\{D_1, \{\!\{ D_2, D_3\}\!\} \}\! \}_L +
\tau_{(123)}\{\!\{D_2, \{\!\{D_3, D_1 \}\!\}\}\!\}_L+
\tau_{(123)}^2 \{\!\{D_3, \{\!\{D_1, D_2\}\!\}\}\}_L=0$$
where 
$$\begin {array}{l}
\{\!\{D_1, \{\!\{ D_2, D_3\}\!\} \}\! \}_L:=
\{\!\{D_1, \{\!\{ D_2, D_3\}\!\}^\prime  \}\! \} \otimes \{\!\{D_2, D_3\}\!\}^{\prime \prime}\\
\forall (D,a) \in \mathbb{L}\times A, \quad \{\!\{D,a \}\!\}=\omega (D)(a)
\end{array}$$
\item If $(D, \Delta ) \in \mathbb{L}^2 $,  one has  

$$\begin{array}{l}
\{\!\{ D, a \Delta \}\!\}= D(a)\Delta +a \{\!\{D, \Delta\}\!\}\\
\{\!\{D,  \Delta a \}\!\}=\{\!\{D, \Delta \}\!\}a+ \Delta D(a)
\end{array}$$
where  $D(a)\Delta = D(a)^\prime  \otimes D(a)^{\prime \prime }\Delta\in A \otimes \mathbb{L}$ and 
$\Delta D(a)= \Delta D(a)^\prime \otimes D(a)^{\prime \prime} \in \mathbb{L}\otimes A$.
\end{itemize}
\end{definition}

\begin{remarks} 
\begin{enumerate}
\item When there will be ambiguity, the bracket $\{\!\{-,-\}\!\}$ will be denoted $\{\!\{-,-\}\!\}_{\mathbb{L}}$.

\item The Jacobi identity is equivalent to the following : 
Let  $D_1 , D_2 , D_3$ be  elements of $\mathbb{L}$ :

$$\{\!\{D_1, \{\!\{ D_2, D_3\}\!\} \}\! \}_R +
\tau_{(123)}\{\!\{D_2, \{\!\{D_3, D_1 \}\!\}\}\!\}_R+
\tau_{(123)}^2 \{\!\{D_3, \{\!\{D_1, D_2\}\!\}\}\}_R=0$$
where 
$\begin {array}{l}
\{\!\{D_1, \{\!\{ D_2, D_3\}\!\} \}\! \}_R:=\{\!\{ D_2, D_3\}\!\}^\prime   \otimes
 \{\!\{D_1,   \{\!\{D_2, D_3\}\!\}^{\prime \prime}\}\!\}.\\
\end{array}$
 \item For differential calculus, we will assume that $A$ is a smooth $B$-algebra and that $\mathbb{L}$ is a finitely generated projective $A^{e}$-module.
 \end{enumerate}
\end{remarks}

{\bf Notation :} We set $\{\!\{ D, \Delta \}\!\}_l$ for the component of $\{\!\{ D, \Delta \}\!\}$ that is 
$\mathbb{L} \otimes A$ and 
$\{\!\{ D, \Delta \}\!\}_r$ for its component in $A \otimes \mathbb{L}$. Adopting a Sweedler's type notation, we set 
$$\begin{array}{l}
\{\!\{ D, \Delta \}\!\}_l= \{\!\{ D, \Delta \}\!\}_l^\prime \otimes \{\!\{ D, \Delta \}\!\}_l^{\prime \prime}\\
 \{\!\{ D, \Delta \}\!\}_r= \{\!\{ D, \Delta \}\!\}_r^\prime \otimes \{\!\{ D, \Delta \}\!\}_r^{\prime \prime}\\
 \end{array}$$
 with $\{\!\{ D, \Delta \}\!\}_l^\prime , \{\!\{ D, \Delta \}\!\}_r^{\prime \prime} \in \mathbb{L}$ and 
 $\{\!\{ D, \Delta \}\!\}_r^\prime , \{\!\{ D, \Delta \}\!\}_l^{\prime \prime} \in A$.  One has 
 $$\{\!\{ D, \Delta \}\!\}_l= -\{\!\{ \Delta , D \}\!\}_r^\circ .$$
 
    \begin{lemma}  \label{Jacobi projected}
 Let $D_1$, $D_2$, $D_3$ be three elements of $\mathbb{L}$. One has : 
$$ \begin{array}{l}
 \{\!\{D_1, \{\!\{D_2, D_3\}\!\}^\prime_l\}\!\}^\prime _l\otimes 
 \{\!\{D_1, \{\!\{D_2, D_3\}\!\}^\prime_l\}\!\}_l^{\prime \prime }
 \otimes \{\!\{D_2, D_3\}\!\}_l^{\prime \prime}\\
 - \{\!\{D_1, D_3\}\!\}^\prime_l \otimes \{\!\{D_2, \{\!\{D_1, D_3\}\!\}^{\prime\prime}_l \}\!\}^\prime \otimes 
 \{\!\{D_2, \{\!\{D_1, D_3\}\!\}^{\prime\prime}_l \}\!\}^{\prime \prime}\\
 -\{\!\{\{\!\{D_1, D_2 \}\!\}^\prime_l, D_3\}\!\}^\prime_l \otimes 
 \{\!\{D_1, D_2\}\!\}_l^{\prime \prime}\otimes
 \{\!\{\{\!\{D_1, D_2 \}\!\}^\prime_l, D_3\}\!\}^{\prime \prime}_l=0
 
 \end{array}$$
 \end{lemma}
 
\begin{proof} The lemma follows from the Jacobi identity, taking the component on 
 $\mathbb{L}\otimes A\otimes A$.\\\end{proof}
 
\begin{proposition}\label{properties of left and right brackets}
Let $( D, \Delta )\in \mathbb{L}^2$ and $(\alpha , \beta ) \in A^2$, one has :\\
$\begin{array}{rcl}
\{\!\{D, \Delta \cdot \beta \}\!\}_l&= &
\{\!\{D, \Delta \}\!\}^\prime_l \otimes \{\!\{D, \Delta \}\!\}^{\prime \prime}_l \cdot \beta + 
\Delta \cdot D(\beta )^\prime \otimes D(\beta)^{\prime \prime}\\
\{\!\{D, \Delta \cdot \beta \}\!\}_r&=&
\{\!\{D, \Delta \}\!\}^\prime _r \otimes \{\{D, \Delta \}\!\}_r^{\prime \prime} \cdot \beta \\
\{\!\{ D, \alpha \cdot \Delta \}\!\}_l&=& \alpha \cdot \{\!\{ D, \Delta \}\!\}_l\\
\{\!\{D, \alpha \cdot \Delta \}\!\}_r&=&
\omega (D)(\alpha)^\prime \otimes \omega (D)(\alpha )^{\prime \prime}\Delta + 
\alpha \{\!\{D, \Delta \}\!\}_r\\
\{\!\{a \cdot D, \Delta\}\!\}_r&=& \{\!\{D, \Delta\}\!\}^\prime _r \otimes a \cdot \{\!\{ D, \Delta\}\!\}^{\prime \prime}_r\\
\{\!\{a \cdot D, \Delta \}\!\}_l&=& 
-\{\!\{\Delta , D \}\!\}^{\prime \prime}_r \otimes a \{\!\{ \Delta , D\}\!\}^\prime _r -
\{\!\{\Delta , a \}\!\}^{\prime \prime}D \otimes \{\!\{\Delta , a \}\!\}^\prime \\
&= &
\{\!\{D, \Delta \}\!\}_l^\prime \otimes a \{\!\{D, \Delta \}\!\}_l ^{\prime \prime}-
\{\!\{\Delta , a \}\!\}^{\prime \prime}D \otimes \{\!\{\Delta , a \}\!\}^\prime \\

\{\!\{D \cdot \beta , \Delta \}\!\}_{r}&= &-\{\!\{\Delta , D \}\!\}^{\prime \prime}_l \beta \otimes \{\!\{ \Delta , D \}\!\}^\prime_l-\{\!\{\Delta , \beta \}\!\}^{\prime \prime}\otimes D \{\!\{\Delta , \beta \}\!\}^\prime\\
\{\!\{D \cdot \beta , \Delta \}\!\}_l&=& -\{\!\{\Delta , D\}\!\}^{\prime \prime}_r \beta \otimes \{\!\{\Delta , D \}\!\}^\prime_r =
 \{\!\{D , \Delta\}\!\}^{\prime }_l \beta \otimes \{\!\{D , \Delta \}\!\}^{\prime \prime}_l

\end{array}$
\end{proposition}

\begin{proof}The proof follows from the properties of the double Lie algebroid  bracket and easy computations.
\end{proof} 
\begin{proposition} \label{characterization Gerstenhaber}
Let $(\mathbb{L}, A, \omega, \{\!\{, \}\!\}_{\mathbb{L}})$ be a quadruple such that 
$\{\!\{,\}\!\}_{\mathbb{L}}$ is a map from $\mathbb{L} \otimes \mathbb{L}$  to
 $A \otimes \mathbb{L} \oplus \mathbb{L} \otimes A$ 
and $\omega$ a map from $\mathbb{L}$ to 
$\mathbb{D}er (A)$.  
There is a unique  graded bilinear map 
$\{\!\{, \}\!\} : T_A(\mathbb{L}) \otimes T_A(\mathbb{L}) \to T_A(\mathbb{L})$ 
of degree $-1$    satisfying the  following conditions : 
\begin{itemize} 
\item For all $(\alpha, \beta, D, \Delta ) \in A^2 \otimes \mathbb{L}^2$, 
$$\begin{array}{l}
\{\!\{\alpha,\beta\}\!\}=0\\
\{\!\{D, \alpha\}\!\}=\omega ( D )(\alpha)^\prime \otimes \omega (D)(\alpha)^{\prime \prime}\\
\{\!\{D, \Delta \}\!\}=\{\!\{D, \Delta \}\!\}_{\mathbb{L}}
\end{array}$$

\item for all $(a,b,c) \in T_A(\mathbb{L})^3$
$$\begin{array}{l}
\{\!\{a, bc\}\!\}= (-1)^{(\vert a \vert -1) \vert b \vert } b\{\!\{a,c\}\!\} + \{\!\{a,b\}\!\}c\\
\{\!\{a,b\}\!\}=-(-1)^{(\vert a \vert -1) (\vert b\vert  -1)}\sigma_{(12)}\{\!\{b,a\}\!\}
\end{array}$$
\end{itemize}

$(\mathbb{L},A, \omega )$ is a double Lie--Rinehart algebra over $A$  if and only if 
$T_A(\mathbb{L})$ is a double Gerstenhaber algebra.
\end{proposition}

\begin{proof}

Adopting a Sweedler type notation for $\{\!\{ , \}\!\}$, one has for any $a_1,\dots ,a_m,b_1,\dots, b_n$ in 
$T_A(\mathbb{L})$~:

$$\{\!\{a_1 \dots a_m, b_1\dots b_n\}\!\} =
{\displaystyle \sum_{p,q}}b_1\dots b_{q-1}\{\!\{a_p, b_q\}\!\}^\prime a_{p+1} \dots a_m \otimes 
a_1\dots a_{p-1}\{\!\{a_p, b_q\}\!\}^{\prime \prime} b_{q+1}\dots b_n
$$

One sets $$\{\!\{a,b,c\}\!\} =\{\!\{a,\{\!\{b,c\}\!\} \}\!\}_{\color{black}L} +
(-1)^{(\vert a \vert -1) (\vert b\vert + \vert c\vert )}
{\color{black}\sigma}_{(123)}\{\!\{b,\{\!\{c,a\}\!\} \}\!\}_{\color{black}L}+ 
(-1)^{\vert c \vert -1) (\vert a \vert + \vert b\vert )}
{\color{black}\sigma}_{(123)}^2\{\!\{c,\{\!\{a,b\}\!\} \}\!\}_{\color{black}L} .$$
$\{\!\{-,-,-\}\!\}$ is a triple bracket. 
We want to show that it is zero. We need to check it on generators. \\

{\it First case :} If two of the $a,b,c$ are in $A$ then all the terms are zero. 

{\it Second case :} If all the $a$, $b$, $c$ are in $\mathbb{L}$, it is $0$ by hypothesis. 

{\it Third case :} We check Jacobi identity in the case of a  triple 
$(\alpha,D,\Delta)$ in $A \times \mathbb{L}\times \mathbb{L}$. 
We need to show the equality : 

$$\{\!\{\alpha,\{\!\{D,\Delta\}\!\} \}\!\}_L +
\sigma_{(123)}\{\!\{D,\{\!\{\Delta,\alpha\}\!\} \}\!\}_L+ 
\sigma_{(123)}^2\{\!\{\Delta,\{\!\{\alpha,D\}\!\} \}\!\}_L =0$$

Before computing each of these terms, let us remark that 
$$\omega \left (\{\!\{ D, \Delta \}\!\}_l \right ) = \{\!\{ \omega (D), \omega (\Delta)\}\!\}_l.$$
This equality can be written more simply as follows :
$$\omega \left ( \{\!\{D, \Delta \}\!\}_l^\prime \right )\otimes \{\!\{D, \Delta\}\!\}^{\prime \prime}_l=
 \{\!\{\omega (D), \omega (\Delta ) \}\!\}_l^\prime \otimes 
 \{\!\{\omega (D), \omega (\Delta)\}\!\}^{\prime \prime}_l.$$
one has 
$$\begin{array}{rcl}
\{\!\{\alpha, \{\!\{ D, \Delta \}\!\} \}\!\}_L&=& 
\{\!\{\alpha, \{\!\{ D, \Delta \}\!\}_l^\prime \}\!\} \otimes \{\!\{ D, \Delta \}\!\}_l^{\prime \prime } \\
&=& -\tau_{(12)}  \omega \left ( \{\!\{D, \Delta \}\!\}^\prime_l \right )(\alpha) \otimes 
\{\!\{D, \Delta \}\!\}^{\prime \prime}_l\\
&=& -\tau_{(12)}   \{\!\{\omega (D), \omega (\Delta ) \}\!\}^\prime_l (\alpha) \otimes 
\{\!\{\omega (D), \omega (\Delta ) \}\!\}^{\prime \prime}_l\\
&=&-\tau_{(123)} \left [ (\omega (D)\otimes 1)\omega (\Delta ) - 
(1 \otimes \omega (\Delta ))\omega (D)\right ] (\alpha)
\end{array}$$
$$\begin{array}{l}
\sigma_{(123)} \{\!\{D, \{\!\{ \Delta ,\alpha \}\!\} \}\!\}_L 
= 
\sigma_{(123)} \left ( \omega (D)
 \left [\omega (\Delta )(\alpha)^\prime \right ] \otimes \omega (\Delta )(\alpha)^{\prime \prime}\right )
= \tau_{(123)} \left ( (\omega (D)\otimes 1)\omega (\Delta) \right )(\alpha)
\end{array}$$
$$\begin{array}{rcl}
\sigma_{(23)} \{\!\{ \Delta , \{\!\{\alpha, D\}\!\} \}\!\}_L
&=& -\tau_{(132)} \left [ \Delta \left ( \omega(D)(\alpha)^{\prime \prime} \right )^\prime 
 \otimes \Delta \left ( \omega(D)(\alpha)^{\prime \prime} \right )^{\prime \prime }
 \otimes \omega (D)(\alpha)^\prime
 \right ]\\
 &=& -\tau_{(123)} \left [\left (1 \otimes \omega (\Delta ) \right ) \omega (D)(\alpha)\right ] 

\end{array}$$\end{proof}
In the next proposition, we study the case where the anchor $\omega$ is zero. 

\begin{proposition}
Let $(\mathbb{L}, \{\!\{-,-\}\!\}, \omega )$ be a double Lie--Rinehart algebra over $A$ such that $\omega=0$. 
Set $X\bullet_rY=-\{X,Y\}_r=-\{\!\{X,Y\}\!\}_r^\prime \{\!\{X,Y\}\!\}_r^{\prime \prime}$ and 
$X\bullet_lY=\{X,Y\}_l=\{\!\{X,Y\}\!\}_l^\prime \{\!\{X,Y\}\!\}_l^{\prime \prime}$.
The laws $\bullet_l$ and $\bullet_r$ are associative so that 
$\{-,-\}:{\mathbb L}\times \mathbb{L} \to {\mathbb{L}}$ is the difference of two associative products. 
The  laws induces by  $\bullet_l$ and $\bullet_r$ in $\dfrac{\mathbb L}{[A, {\mathbb L}]}$
are opposite from each other. The bracket $\{-,-\}$ induces a Lie bracket on 
$\dfrac{\mathbb L}{[A, {\mathbb L}]}$ that  comes from an associative product. 
\end{proposition}
\begin{proof}
The associativity of $\bullet_l$ follows from lemma \ref{Jacobi projected} and 
proposition \ref{properties of left and right brackets}. The proof of the associativity of $\bullet_r$ is similar. 
Let $X$ and $Y$ in $\mathbb L$. 
$$\begin{array}{rcl}
X\bullet_r Y&=&- \{\!\{X,Y\}\!\}_r^\prime \{\!\{X,Y\}\!\}_r^{\prime \prime}
= \{\!\{Y,X\}\!\}_l^{\prime \prime} \{\!\{Y,X\}\!\}_l^{\prime }= 
 \{\!\{Y,X\}\!\}_l^{\prime } \{\!\{Y,X\}\!\}_l^{\prime \prime}\; {\rm mod}\; [A,{\mathbb L}]\\
 &=&Y\bullet_l X\; {\rm mod}\; [A,{\mathbb L}].
 \end{array}$$
\end{proof}
 
 Let us now give several examples of double Lie algebroids. \\

{\bf Example 1 :} 

$(\mathbb{D}er (A) , A, id)$ is a double Lie  algebroid. \\

{\bf Example 2 :} If $A=k$, double Lie algebroid structures on $\mathbb{L}$ over $k$ are in bijection with associative $k$-algebras structures over $\mathbb{L}$. Indeed, If $\mathbb{L}$ is a double Lie  algebroid over $k$, then 
$\{\!\{-,-\}\!\}_l$ is an associative product. Conversely, an associative product over $k$ gives rise to a double Lie--Rinehart algebra over $k$ as follows :
$\{\!\{X,Y\}\!\}_l=XY\otimes 1$, $\{\!\{X,Y\}\!\}_r=-1\otimes YX$ and  $\{\!\{X,Y\}\!\}=XY\otimes 1 -1 \otimes YX$. The bracket $\{\!\{-,-\}\!\}$ identifies to the Lie bracket coming from the associative structure. \\


{\bf Exemple 3 :}

Let us recall the definition of a double Lie  algebra :

\begin{definition} Let $\mathfrak g$ be a vector space of finite dimension. A double Lie bracket over 
$\mathfrak g$ is a map $\{\!\{-,-\}\!\} : {\mathfrak g} \otimes {\mathfrak g} \to {\mathfrak g} \otimes {\mathfrak g}$ 
satisfying the Jacobi identity.
\end{definition}

\begin{remark} Let $V$ be a $m$ dimensional $\mathbb{C}$-vector space. Double Lie algebra structures on $V$ are in bijection  (\cite{A-R-R}) with  operators  $r\in End(V\otimes V)$ such that 
$r(v\otimes u)=-r( u \otimes v)^\circ$ (skew-summetry)  and solution of the classical Yang Baxter  equation 
$$r^{23}r^{12}+r^{31}r^{23}+r^{12}r^{31}=0$$
where $r^{ij}$ acts on $V^{\otimes 3}$ non trivially on $(i,j)$ spaces and as identity elsewhere. \\

\end{remark}

Let $\mathfrak g$ be a double Lie algebra with double bracket 
$\{\!\{-,-\}\!\}_{\mathfrak g}$. This latter induces a Poisson bracket on $T({\mathfrak g})$. There is a unique double Lie algebroid structure on 
$\mathbb{L}=T({\mathfrak g})\otimes {\mathfrak g}\otimes T({\mathfrak g})$ such that 
$$\begin{array}{l} 
\forall (X,Y)\in {\mathfrak g}, \quad \forall a \in T({\mathfrak g}), \quad \\
\{\!\{X,Y\}\!\}_{\mathbb{L}}=\{\!\{X,Y\}\!\}_{\mathfrak g}\\
\omega (X)(a)=\{\!\{X,a\}\!\}_{T({\mathfrak g})}
\end{array}$$
Then $T({\mathfrak g}) \otimes {\mathfrak g} \otimes T({\mathfrak g})$ is a Lie double algebroid over 
$T({\mathfrak g})$. \\

{\bf Example 4:} Let $V$ be an $A^e$-module of finite type. For 
$\lambda$, $\mu$ in $Hom_{A^e}(V, V \otimes V)$ (the exterior $A^e$-module structure on $V \otimes V$ is used for the $Hom_{A^e}$), one sets :
$$\begin{array}{l}
\{\!\{\lambda , \mu \}\!\}_{l}^{\sim}= (\lambda \otimes 1)\mu -(1 \otimes \mu )\lambda \\
\{\!\{\lambda , \mu \}\!\}_r^{\sim}= (1 \otimes \lambda)\mu -(\mu \otimes 1 )\lambda \\
\end{array}$$
One sets 
$$\begin{array}{l}
\{\!\{\lambda , \mu \}\!\}_l= \tau_{(23)}\{\!\{\lambda , \mu \}\!\}^{\sim}_l\\
\{\!\{\lambda , \mu \}\!\}_r= \tau_{(12)}\{\!\{\lambda , \mu \}\!\} ^{\sim}_r\\
\{\!\{\lambda , \mu \}\!\}=\{\!\{\lambda , \mu \}\!\}_l + \{\!\{\lambda , \mu \}\!\}_r
\end{array}$$
Then 
$$\begin{array}{rcl}
\{\!\{\lambda , \mu \}\!\} \in Hom_{A^e}(V, V \otimes V) \otimes_k V + 
V \otimes _kHom_{A^e}(V, V \otimes V)\\
\end{array}$$
Then 
$$\begin{array}{rcl}
\{\!\{\lambda , \mu \}\!\} \in Hom_{A^e}(V, V \otimes V) \otimes_k T_{A}(V) \oplus 
T_A(V)\otimes_k Hom_{A ^e}(V, V \otimes V)
\end{array}$$
There is a unique double Lie algebroid structure on 
$T_A(V) \otimes_{k} Hom (V, V \otimes V) \otimes_{k} T_A(V)$
with anchor map 
$$\begin{array}{rcl}
\omega :  {\color{black}T_A(V) \otimes_{k} Hom (V, V \otimes V) \otimes_{k} T_A(V)}
 &\to &\mathbb{D}er (T(V))\\
\alpha \otimes \lambda \otimes \beta & \mapsto & \alpha\{\!\{\lambda , -\}\!\} \beta
\end{array}$$
where $\{\!\{\lambda , -\}\!\} $ is the unique {\color{black} double} 
derivation of $T_A(V)$ such that for all $v\in V$, 
$\{\!\{\lambda, v\}\!\}=\lambda (v)$ and, for all $a\in A$, $\{\!\{\lambda, a\}\!\}=0$.\\

{\bf Example 5 :}

Let $A$ be a $B$-algebra endowed with a $B$-linear double Poisson bracket.  
The quadruple $\left ( \Omega_B^1A , A, \omega , \{\!\{ , \}\!\} \right )$  is a double Lie--Rinehart algebra with 
$\omega$ and $\{\!\{, \}\!\}$ defined as follows :
 
$$\begin{array}{rcl}
\Omega_B^1A& \to & {\mathbb D}er_B (A)\\
da& \mapsto & \{\!\{a,- \}\!\}
\end{array}$$
$$\{\!\{da, db\}\!\}=d \left (\{\!\{a,b\}\!\}^{\prime} \right ) \otimes \{\!\{a,b\}\!\}^{\prime \prime} +
\{\!\{a,b\}\!\}^\prime \otimes d \left (\{\!\{a,b\}\!\}^{\prime\prime} \right )$$\\

{\bf Example 6 :} (\cite{VdB2}) Let $(A,P)$ a double quasi-Poisson algebra
{\color{black} and let $E$ be the double derivation of $A$  defined by $E(a)=a\otimes 1-1 \otimes a$}
 Then 
$\tilde{\Omega}_A=\Omega_A \oplus AEA$ has the structure of a double Lie algebroid where the double bracket is defined as follows
$$\begin{array}{l}
\{\!\{da,b\}\!\}_{\tilde{\Omega}A}=\{\!\{a,b\}\!\}\\
\{\!\{da,db\}\!\}_{\tilde{\Omega}A}=d\{\!\{a,b\}\!\}+ 
\dfrac{1}{4}\left [ b, \left [a, E \otimes 1 -1 \otimes E\right ]_*\right ]\\
\{\!\{E,X\}\!\}_{\tilde{\Omega}A}=X \otimes 1-1 \otimes X
\end{array}$$
for $a,b\in A, \quad X \in T_A \tilde{\Omega}_A$ where $[-,-]_*$ denotes the commutator for the inner $A$-bimodule structure on $AEA \otimes AEA$. Futhermore the anchor is the A$^e$-bimodule morphism 
defined by : 
$$\begin{array}{rcl}
\Omega_A \oplus AEA &\to& \mathbb{D}er(A)\\
(du , \delta ) & \mapsto & \{\!\{u,-\}\!\} +\delta\quad 
\end{array}$$
is surjective.

\section{Differential calculus}



$A \otimes A$ is endowed with a $A^e \otimes A^e$-module structure. 
Let $M$ be an $A^e$ -module. We choose to  set 
$$M^*=\{\lambda :M \to A\otimes A \mid 
\lambda (\alpha \cdot D \cdot \beta )= \lambda (D)^\prime \beta  \otimes \alpha \lambda (D)^{\prime \prime} \}$$
$M^*$ is itself an $A^e$-module as follows : 
$$\forall \lambda \in M^*, \quad \forall (a,b)\in A^2, \quad \forall D \in M, \quad 
(a \cdot \lambda \cdot b)(D)=a \lambda (D)^\prime \otimes \lambda (D)^{\prime \prime}b.$$

\begin{remark}
We can  exchange the role of the two $A^e$-module structures on $A \otimes A$ and define 
$$M_*= \{\lambda : M \to A \otimes A \mid 
\lambda (\alpha \cdot D \cdot \beta )= 
\alpha  \lambda (D)^\prime  \otimes  \lambda (D)^{\prime \prime}\beta \}.$$
Composition by $\tau_{(12)}$ is an isomorphism of $A^e$-modules from $M^*$ to $M_*$. \end{remark}

Let  $M$ is an $A^e$- module. Let us endow $A^{\otimes n+1}$ with the $(A^e)^{\otimes n}$-module structure where the ith copy of $A^e$ acts as follows : 
$$a^\prime \otimes a^{\prime \prime} \cdot (a_1 \otimes \dots \otimes a_{n+1})=
a_{1}\otimes \dots \otimes a_i a^{\prime \prime}\otimes a^\prime a_{i+1}\otimes \dots \otimes a_{n+1},$$
then the map (\cite{VdB1})
\footnote{our map is slightly different from that of M. Van den Bergh  due to different conventions for the dual of 
$\mathbb{L}$ : 
\cite{VdB1} makes use of $M_*$, we make use of $M^*$.} 
$$\Psi : M^{*\otimes _A n}  \to  Hom_{(A^e)^{\otimes n}} \left (M^{\otimes n}, A ^{\otimes n+1} \right )$$

$$\Psi (\lambda_1 \otimes \lambda_2 \otimes \dots \otimes \lambda_n) 
(m_1 \otimes \dots \otimes m_n )=
\lambda_1(m_1)^\prime \otimes \lambda _1(m_1)^{\prime \prime}\lambda_2 (m_2 )^\prime \otimes \dots \otimes 
\lambda_{n-1} (m_{n-1} )^{\prime \prime}\lambda_n (m_n)^\prime \otimes \lambda_n (m_n)^{\prime \prime}
$$
is well defined. 

If $M$ is a finitely generated $A^e$-projective module, then $\Psi$ is an isomorphism of $A^e$-modules. \\

The cyclic group $C_n$ acts on $Hom_{(A^ e)^{\otimes n}} \left ( M ^{ \otimes n},  A^{\otimes n}\right )$ as follows :
$$\forall \omega \in Hom_{(A^ e)^{\otimes n}} \left ( M ^{ \otimes n},  A^{\otimes n}\right ), \quad 
\tau_{(1\dots n)}\cdot \omega = \tau_{(1\dots n )}\circ \omega \circ \tau_{(1\dots n)}^{-1}$$
The set of signed invariants of $Hom_{(A^ e)^{\otimes n}} \left ( M ^{ \otimes n},  A^{\otimes n}\right )$ under the action of $C_n$ is
$$s-inv \;Hom_{(A^ e)^{\otimes n}} \left ( M ^{ \otimes n},  A^{\otimes n}\right ):=
\{ \omega \in Hom_{(A^ e)^{\otimes n}} \left ( M ^{ \otimes n},  A^{\otimes n}\right ) \mid 
\tau_{(1\dots n )} \cdot \omega = (-1)^{n-1}\omega \}$$

If $M$ is a finitely generated projective $A^e$-module,  M. van den Bergh showed (\cite{VdB1})  that 
$\dfrac{T^*(M^*)}{[T^*(M^*), T^*(M^*)]}$ is isomorphic to 
$s-inv  Hom_{(A^ e)^{\otimes n}} \left ( M ^{ \otimes n},  A^{\otimes n}\right )$. He constructed  the following  isomorphism $\mu$ between the two spaces :

$\mu : \dfrac{T^*(M^*)}{[T^*(M^*), T^*(M^*)]} \to  
s-inv \; Hom_{(A^ e)^{\otimes n}}  \left ( M ^{ \otimes n},  A^{\otimes n}\right )$
$$\mu (\lambda_1 \otimes \dots \otimes \lambda_n)=
\{\!\{-, \dots , -\}\!\}_{\lambda_1 \dots \lambda_n}=\sum_{i}(-1)^{(n-1)i} \tau_{(1\dots n )}^i \circ  \Phi \circ  \tau^{-i}_{(1\dots n)}$$ 
with  
$$\Phi (\lambda_1 \otimes \lambda_2 \otimes \dots \otimes \lambda_n) (m_1 \otimes \dots \otimes m_n )=
\lambda_n (m_n)^{\prime \prime}\lambda_1(m_1)^\prime \otimes
 \lambda _1(m_1)^{\prime \prime}\lambda_2 (m_2 )^\prime \otimes \dots \otimes 
\lambda_{n-1} (m_{n-1} )^{\prime \prime}\lambda_n (m_n)^\prime .$$
\vspace{1em}

{\bf Notation :}  In the computation, 
$\Phi \left ( \lambda_1 \otimes \dots \otimes \lambda_n \right )$
 will be denoted
 $\widetilde{\lambda_1 \otimes \dots \otimes \lambda_n}$.\\
 
Let $C$ be a graded algebra  and let 
$$\begin{array}{rcl}
\Theta : C & \to & C \otimes C \\
c & \mapsto& \Theta (c)^\prime \otimes \Theta (c)^{\prime \prime}
\end{array}$$
be a double derivation. One sets 
$$\begin{array}{rcl}
^\circ  \Theta : C & \to & C \\
c & \mapsto& \Theta (c)^{\prime \prime } \Theta (c)^{\prime }.
\end{array}$$
$^\circ \Theta$ is an endomorphism of $C$ and induces an endomorphism of  $\dfrac{C}{[C, C]}$. 
We {\color{black} will} be mostly interested in the case where $C=T_A(\mathbb{L}^*)$ or 
$C=T_A(\mathbb{L})$ 
where $\mathbb{L}$ is a double Lie--Rinehart algebra.\\

\begin{proposition} \label{reduced and unreduced}
1) Let $d: C \to C$ be a derivation of degree $\mid d \mid $. It induces a derivation still denoted 
$d: C \otimes C \to C\otimes C$ :
$$\forall c_1\otimes c_2 \in C \otimes C, \quad 
d(c_1 \otimes c_2)=d(c_1)\otimes c_2 + (-1)^{\mid c_1 \mid} c_1 \otimes d(c_2)$$
Let  $\underline{i} : C \to C \otimes C$ be a double derivation.  
If  $c\in C$, we will write $\underline{i}(c)= \underline{i} (c)^\prime \otimes 
\underline{i}(c)^{\prime \prime}$ and we set $\iota={}^\circ\underline{i}$. 
Then ${}^\circ(d\circ \underline{i}) = d \circ \iota$ and 
${}^\circ ( \underline{i} \circ d ) = \iota \circ d$.

2)( \cite{VdB1}) Let $\delta : C \to C \otimes C$ and $\Delta :C \to C\otimes C$ be two double derivations. 
$$\delta \circ ^\circ \Delta -\tau_{(12)} \Delta \circ ^\circ \delta=^{\circ ,l}\{\!\{\delta , \Delta\}\!\}_l +
^{\circ ,r}\{\!\{\delta , \Delta\}\!\}_r$$
where 
$^{\circ ,l}(\epsilon^\prime \otimes \epsilon^{\prime \prime})=^\circ \epsilon^\prime
\otimes \epsilon^{\prime \prime}$ and 
$^{\circ ,r}(\epsilon^\prime \otimes \epsilon^{\prime \prime})= \epsilon^\prime
\otimes ^\circ \epsilon^{\prime \prime}$
\end{proposition}

\begin{proof}
1) $(d \circ \underline{i})(c)= d \left ( \underline{i}(c)^{\prime} \right ) \otimes 
\underline{i}(c)^{\prime \prime} + 
 (-1)^{\left \vert  \underline{i}(c)^{\prime}\right \vert} 
  \underline{i}(c)^{\prime} \otimes d \left (  \underline{i}(c)^{\prime \prime}\right )$.\\
Hence 
$$\begin{array}{rcl}
{}^\circ(d \circ \underline{i}) (c)&=&
 \underline{i}(c)^{\prime \prime} d \underline{i}(c)^{\prime}
 (-1)^{(\vert  \underline{i}(c)^{\prime}\vert +\mid d \mid )\vert  \underline{i}(c)^{\prime \prime} \vert}
 +\left ( d \underline{i}(c)^{\prime \prime} \right ) \underline{i}(c)^{\prime}
 (-1)^{\left \vert  \underline{i}(c)^{\prime \prime}\right \vert \left \vert  \underline{i}(c)^{\prime}\right \vert }\\
 &=& d \left [
 \underline{i}(c)^{\prime \prime}\underline{i}(c)^{\prime}
 (-1)^{\vert \underline{i}(c)^{\prime}\vert \vert \underline{i}(c)^{\prime \prime}\vert }
  \right ]\\
  &=& (d \circ \iota )(c)
\end{array}$$

The equality  ${^\circ}(\underline{i} \circ d ) = \iota \circ d$ is obvious. 

2) is stated in \cite{VdB1}.
\end{proof}

Let us now see examples of this situation. \\

{\bf The contraction } (similar to \cite{CB-PE-VG})

Set $\mathbb{L}^*=Hom_{A^e}(\mathbb{L}, A \otimes A)$ and 
let  $D \in \mathbb{L}$. The element  $D$ defines a degree {\color{black}-1} double  derivation 
$\underline{i}_{D}: T_A(\mathbb{L}^*) \to T_A(\mathbb{L}^*)\otimes T_A(\mathbb{L}^*)$
$$\begin{array}{l}
\forall \alpha \in \mathbb{L}^*, \quad \underline{i}_D (\alpha)=\alpha (D)^\prime \otimes  \alpha (D)^{\prime \prime}.
\end{array}$$
More explicitely, the map 
$\underline{i}_D: T^\bullet (\mathbb{L}^*) \to T ^\bullet (\mathbb{L}^*)\otimes T^\bullet (\mathbb{L}^*)$ 
is given by
$$\underline{i}_D(\alpha_1 \dots \alpha_n)=
\displaystyle \sum_{i=1}^n (-1)^{k-1}
\alpha_1 \dots \alpha_{k-1}\alpha_k(D)^\prime \otimes \alpha_k(D)^{\prime \prime} \alpha_{k+1} \dots \alpha_n $$
\begin{lemma}(\cite{CB-PE-VG})
For all $\Phi , \Theta$ in $\mathbb{L}$, one has :
$$\begin{array}{l}
\underline{i}_\Theta (\alpha \beta)=\underline{i}_{\Theta}(\alpha ) \beta +
(-1)^{deg \; \alpha} \alpha \underline{i}_{\Theta}(\beta )\\
\underline{i}_\Phi  \circ \underline{i}_{\Theta}+ \underline{i}_\Theta  \circ \underline{i}_{\Phi}=0
\end{array}$$
\end{lemma}
\vspace{1em}

{\it From now on, we will assume that   ${\mathbb L}$ is a finitely generated and projective $A^{e}$-module}\\

 \vspace{1em}

The  Karoubi-de Rham complex was defined in \cite{Karoubi}. It corresponds to the case 
where  $\mathbb{L}=\mathbb{D}er(A)$. In this case, 
differential calculus was treated in \cite{CB-PE-VG} (differential, Lie derivative etc...) but   the formulas don't adjust to any double Lie algebroid. Even if  $\mathbb{L}=\mathbb{D}er(A)$, some of our formulas are different but, most of the time, we make use of the hypothesis $A$ smooth. \\

{\bf The differential }

For  a general double Lie algebroid, the definition of the differential is more complicated than in the case where $\mathbb{L}=\mathbb{D}er (A)$. 

\begin{theorem}
One defines $\underline{d}_\mathbb{L}: T_A(\mathbb{L}^*) \to T_A(\mathbb{L}^*)$ as the degree one derivation determined by :  
\begin{itemize}
\item $ \forall D \in \mathbb{L}, \quad \underline{d}_\mathbb{L} a (D)=D (a)^\prime \otimes D (a)^{\prime \prime}$
\item $\forall \lambda \in \mathbb{L}^* , \quad 
\underline{d}_{\mathbb{L}}\lambda \in \mathbb{L}^* \otimes _A \mathbb{L}^*\simeq 
Hom_{A^e} (\mathbb{L}\otimes \mathbb{L}, A^{\otimes 3})$ is given by for all : $ D, \Delta  \in {\mathbb L}$
$$\begin{array}{rcl}
\underline{d}_\mathbb{L} (\lambda )(D,\Delta)&=& 
D\left ( \lambda (\Delta) ^\prime\right ) \otimes \lambda (\Delta )^{\prime \prime} -
\lambda (D)^\prime \otimes \Delta \left ( \lambda (D)^{\prime \prime}\right ) -
\tau_{(23)}\left [\lambda (\{\!\{D,\Delta \}\!\}^\prime_l) \otimes \{\!\{D, \Delta \}\!\}^{\prime\prime}_l\right ]\\
\end{array}$$
One can easily compute the general formula for $\underline{d}_{\mathbb{L}}$ : 
$\forall \psi \in Hom_{A^e} (\mathbb{L}^{\otimes n}, A^{\otimes n+1})$
\end{itemize}
$$\begin{array}{rcl}
 \underline{d}_\mathbb{L}(\psi) (D_1, \dots , D_{n+1})&=&
 {\displaystyle \sum_{i=1}^{n+1}}(-1)^{i-1}
 \left ( id^{i-1}\circ D_i \circ id^{n-i}\right ) \psi (D_1, \dots , \widehat{D_i} \dots D_{n+1})\\
 &+& {\displaystyle \sum_{i=1}^n}(-1)^{i}
 \psi \left ( D_1\otimes \dots \otimes \{\!\{ D_i , D_{i+1}\}\!\}_l \otimes \dots \otimes D_{n+1}\right )
\end{array}$$
where if $\psi =\psi_1 \otimes \psi_2 \otimes \dots  \otimes \psi_n
{\color{black}\in{\mathbb{L}^{*\otimes_A n}}}$ with
$\psi _i: \mathbb{L} \to A^{\otimes 2}$, then
$$\psi_{i}(\{\!\{D_i, D_{i+1}\}\!\}_l)=\tau_{(23)} \left [ 
\psi_i (\{\!\{D_i, D_{i+1}\}\!\}_l^\prime) \otimes 
\{\!\{D_i, D_{i+1}\}\!\}_l^{\prime \prime }\right ]$$

  $\underline{d}_{\mathbb{L}^*}$ is a derivation whose  square is zero 
 (that is to say a differential).
 
  \end{theorem}
 

  
 
 

  {\it Proof :} In the proof, we will make use of the following notation : If $\lambda \in \mathbb{L}^*$, then 
  $\underline{d}_\mathbb{L}(\lambda)= 
  \underline{d}_\mathbb{L}(\lambda)^{(1)}\otimes \dots \otimes \underline{d}_\mathbb{L}(\lambda)^{(n+1)}$ 
  where 
  $\underline{d}_\mathbb{L}(\lambda)^{(i)}$ takes values in $A$. \\
  
  If $D_1$, $D_2$, $D_3$ are three elements of $\mathbb{L}$, one has :
  
 $$\begin{array}{rcl}
( \underline{d}_\mathbb{L} \circ \underline{d}_\mathbb{L})(\lambda)(D_1, D_2, D_3)
&=& 
 D_1\left [ D_2 \left ( \lambda (D_3)^\prime \right )^\prime \right ] \otimes 
 D_2 \left ( \lambda (D_3)^\prime \right )^{\prime \prime} \otimes \lambda (D_3)^{\prime \prime} 
 -D_1 \left ( \lambda (D_2)^\prime \right ) \otimes D_3 \left ( \lambda (D_2)^{\prime \prime} \right ) \\
&- &
D_1 \left [ \lambda \left ( \{\!\{D_2, D_3\}\!\}_l^\prime\right )^\prime \right ] 
\otimes \{\!\{D_2, D_3\}\!\}_l^{\prime \prime} \otimes 
\lambda \left ( \{\!\{D_2, D_3\}\!\}_l^\prime\right )^{\prime \prime}\\
&-& 
{\color{black}D_1 \left ( \lambda (D_3)^\prime \right )^\prime \otimes 
D_2 \left [ D_1 \left ( \lambda (D_3)^\prime \right )^{\prime \prime}\right ] \otimes 
 \lambda (D_3)^{\prime \prime}}\\
&+& \lambda (D_1)^\prime \otimes 
D_2 \left [ D_3 \left ( \lambda (D_1)^{\prime \prime}\right )^\prime\right ] \otimes 
D_3 \left ( \lambda (D_1)^{\prime \prime}\right )^{\prime \prime}\\
&+& \lambda \left ( \{\!\{D_1, D_3\}\!\}_l^\prime\right )^\prime 
\otimes D_2 \left [ \{\!\{D_1, D_3\}\!\}_l^{\prime \prime}\right ] \otimes 
\lambda \left ( \{\!\{D_1, D_3\}\!\}_l^\prime\right )^{\prime \prime}\\
 &+&D_1 \left ( \lambda (D_2)^\prime \right ) \otimes D_3 \left ( \lambda (D_2)^{\prime \prime} \right )\\
 &-& \lambda (D_1)^\prime \otimes 
 D_2 \left ( \lambda (D_1)^{\prime \prime}\right )^\prime \otimes 
D_3 \left [ D_2 \left ( \lambda (D_1)^{\prime \prime}\right )^{\prime \prime}\right ]\\
&-& \lambda \left ( \{\!\{D_1, D_2\}\!\}_l^\prime\right )^\prime 
\otimes  \{\!\{D_1, D_2\}\!\}_l^{\prime \prime} \otimes 
D_3 \left [\lambda \left ( \{\!\{D_1, D_2\}\!\}_l^\prime\right )^{\prime \prime}\right ]\\
&-& (\underline{d}_\mathbb{L}\lambda)^{(1)}\left ( \{\!\{D_1, D_2\}\!\}^\prime_l , D_3 \right ) \otimes 
\{\!\{D_1, D_2\}\!\}^{\prime\prime}_l\otimes (\underline{d}_\mathbb{L}\lambda)^{(2)}(-) \otimes 
(\underline{d}_\mathbb{L}\lambda)^{(3)}(-)\\
&+& (\underline{d}_\mathbb{L}\lambda)^{(1)}\left ( D_1, \{\!\{D_2,D_3\}\!\}^\prime_l \right ) \otimes 
 (\underline{d}_\mathbb{L}\lambda)^{(2)} (-)
\otimes \{\!\{D_2, D_3\}\!\}^{\prime\prime}_l \otimes 
(\underline{d}_\mathbb{L}\lambda) ^{(3)}(-)
 \end{array}$$
 But
 $$\begin{array}{l} 
 - (\underline{d}_\mathbb{L}\lambda)^{(1)}\left ( \{\!\{D_1, D_2\}\!\}^\prime_l , D_3 \right ) \otimes 
\{\!\{D_1, D_2\}\!\}_l^{\prime\prime}\otimes (\underline{d}_\mathbb{L}\lambda)^{(3)}(-) \otimes 
(\underline{d}_\mathbb{L}\lambda)^{(4)}(-)\\
+ (\underline{d}_\mathbb{L}\lambda)^{(1)}\left ( D_1, \{\!\{D_2,D_3\}\!\}^\prime_l \right ) \otimes 
 (\underline{d}_\mathbb{L}\lambda)^{(2)} (-)
\otimes \{\!\{D_2, D_3\}\!\}^{\prime\prime}_l \otimes 
(\underline{d}_\mathbb{L}\lambda) ^{(4)}(-)\\
=
{\color{black}- \{\!\{D_1, D_2 \}\!\}^\prime_l\left ( \lambda (D_3)^\prime \right )^\prime \otimes 
\{\!\{D_1,D_2\}\!\}^{\prime \prime}_l \otimes 
 \{\!\{D_1, D_2 \}\!\}^\prime_l\left ( \lambda (D_3)^\prime \right )^{\prime \prime}
\otimes \lambda (D_3)^{\prime \prime}}\\
+\lambda \left (\{\!\{D_1,D_2\}\!\}_l ^\prime \right)^\prime \otimes 
\{\!\{D_1, D_2 \}\!\}^{\prime \prime} \otimes 
D_3 \left [ \lambda \left (\{\!\{D_1,D_2\}\!\}_l ^\prime \right)^\prime\right ]\\
+ \lambda \left[ \{\!\{ \{\!\{D_1,D_2\}\!\}^\prime_l, D_3 \}\!\}_l^\prime \right ]^\prime
\otimes \{\!\{D_1,D_2\}\!\}_l^{\prime \prime}\otimes 
 \{\!\{ \{\!\{D_1,D_2\}\!\}^\prime_l, D_3 \}\!\}_l^{\prime \prime} \otimes 
 \lambda \left [  \{\!\{ \{\!\{D_1,D_2\}\!\}^\prime_l, D_3 \}\!\}_l^\prime\right ]^{\prime \prime}\\
 +D_1 \left [ \lambda \left (\{\!\{D_2,D_3\}\!\}_l ^\prime \right)^\prime \right ] \otimes 
\{\!\{D_2, D_3 \}\!\}_l^{\prime \prime} \otimes 
\lambda \left (\{\!\{D_2,D_3\}\!\}_l ^\prime \right)^{\prime\prime}\\
-\lambda (D_1)^\prime \otimes \{\!\{D_2, D_3\}\!\}_l^\prime \left [ \lambda (D_1)^{\prime \prime}\right ]^\prime
\otimes  \{\!\{D_2, D_3\}\!\}_l^{\prime \prime}
\otimes \{\!\{D_2, D_3\}\!\}_l^\prime \left [ \lambda (D_1)^{\prime \prime}\right ]^{\prime\prime}\\
- \lambda \left [ \{\!\{D_1, \{\!\{D_2, D_3\}\!\}^\prime_l \}\!\}^\prime_l\right ]^\prime
\otimes  \{\!\{D_1, \{\!\{D_2, D_3\}\!\}^{ \prime}_l \}\!\}_l^{\prime \prime} \otimes 
\{\!\{D_2, D_3\}\!\}_l^{\prime \prime}\otimes
 \lambda \left [ \{\!\{D_1, \{\!\{D_2, D_3\}\!\}^\prime_l \}\!\}^\prime_l\right ]^{\prime\prime}
 \end{array}$$
 A lot of terms cancel, and we are left with
 $$\begin{array}{l}
(\underline{d}_\mathbb{L} \circ \underline{d}_\mathbb{L})(\lambda) (D_1, D_2, D_3)=
   \lambda \left ( \{\!\{D_1, D_3\}\!\}_l^\prime\right )^\prime 
\otimes D_2 \left [ \{\!\{D_1, D_3\}\!\}_l^{\prime \prime}\right ] \otimes 
\lambda \left ( \{\!\{D_1, D_3\}\!\}_l^\prime\right )^{\prime \prime}\\
+  \lambda \left[ \{\!\{ \{\!\{D_1,D_2\}\!\}^\prime_l, D_3 \}\!\}_l^\prime \right ]^\prime
\otimes \{\!\{D_1,D_2\}\!\}^{\prime \prime}\otimes 
 \{\!\{ \{\!\{D_1,D_2\}\!\}^\prime_l, D_3 \}\!\}_l^{\prime \prime} \otimes 
 \lambda \left [  \{\!\{ \{\!\{D_1,D_2\}\!\}^\prime_l, D_3 \}\!\}_l^\prime\right ]^{\prime \prime}\\
 - \lambda \left [ \{\!\{D_1, \{\!\{D_2, D_3\}\!\}^\prime_l \}\!\}^\prime_l\right ]^\prime
\otimes  \{\!\{D_1, \{\!\{D_2, D_3\}\!\}^{ \prime}_l \}\!\}_l^{\prime \prime} \otimes 
\{\!\{D_2, D_3\}\!\}_l^{\prime \prime}\otimes
 \lambda \left [ \{\!\{D_1, \{\!\{D_2, D_3\}\!\}^\prime_l \}\!\}^\prime_l\right ]^{\prime\prime}

  \end{array}$$
  
  The equality $\underline{d}_\mathbb{L} \circ \underline{d}_\mathbb{L} (\lambda)(D_1, D_2, D_3)=0$ 
  follows now from the lemma \ref{Jacobi projected}.$\Box$\\

 \begin{remark}
 In \cite{CB-PE-VG}, it is shown that for $\mathbb{L}=\mathbb{D}er (A)$ and $A$ smooth, the complex 
 $\left ( T_A^\bullet (\mathbb{L}), \underline{d}_{\mathbb{L}}\right )$ is acyclic in strictly positive degree.
 \end{remark}

 

 \begin{definition}
 The differential $\underline{d}_\mathbb{L}$ induces a differential 
 $$d_\mathbb{L} : 
 \dfrac{T_A(\mathbb{L}^*)}{[T_A(\mathbb{L}^*), T_A(\mathbb{L}^*)]} \to
  \dfrac{T_A(\mathbb{L}^*)}{[T_A(\mathbb{L}^*),T_A(\mathbb{L}^*)]}.$$
  
 We set 
  $DR(\mathbb{L})^\bullet :=
 \dfrac{T_A(\mathbb{L}^*)}{[T_A(\mathbb{L}^*), T_A(\mathbb{L}^*)]} $ 
 and $d_\mathbb{L}$ will be called 
 "the differential of the double Lie algebroid" $\mathbb{L}$. 
  
  \end{definition}
  
   \begin{remarks}
 
 1) If $\mathbb{L}={\mathbb {D}}er_k (A)$ and $A$ smooth, the complex 
 $\left (DR^\bullet (\mathbb{L}), d_\mathbb{L} \right )$ is 
  the Karoubi--de Rham complex (\cite{Karoubi}, \cite{CB-PE-VG}).
 
 2) If $A$ is a smooth double Poisson algebra and $\mathbb{L}=\Omega^1_k(A)$, we will see that 
the complex 
 $\left (DR^\bullet (\mathbb{L}), d_\mathbb{L} \right )$ is the complex computing   the non commutative Poisson cohomology
 (\cite {Pichereau-vanWeyer}). 
 \end{remarks}
  
 In low degree,  the  expression of $d_\mathbb{L}$ is the following : For all $a\in A$, 
 $\phi \in \mathbb{L}$ and
 $\widetilde{\phi} \in s-inv Hom _{(A^e)} \left ( \mathbb{L} , A \right )$
 $$\begin{array}{l}
 d_\mathbb{L}(a)(D)=D(a)^{\prime \prime} D(a)^\prime \\
 d_\mathbb{L}(\widetilde{\phi} )(D, \Delta )=D\left ( \widetilde{\phi} (\Delta )\right )  
 - \Delta \left (\widetilde{ \phi} (D)\right )^\circ-
 \left [ \widetilde{\phi} \left ( \{\!\{D, \Delta \}\!\}_l^\prime\right ) \otimes \{\!\{D, \Delta \}\!\}_l^{\prime \prime} \right ] 
 -\left [ \{\!\{D, \Delta \}\!\}_r^\prime \otimes \widetilde{\phi} \left ( \{\!\{D, \Delta \}\!\}_r^{\prime \prime} \right )\right ] \\

 \end{array}$$\\
 
 



We now give an expression for $d_L$ in any degree. \\

\begin{theorem}\label{explicit formula}
Let $\mu_1, \dots ,\mu_n \in \mathbb{L}^*$. 
 We have the following formula where indices should be understood modulo $n+1$ : 
$$\begin{array}{l}
d_\mathbb{L} \left (\{\!\{-, \dots , -\}\!\}_{\mu_1 \dots \mu_n} \right )(D_1 \otimes \dots \otimes D_{n+1})\\
={\displaystyle \sum_{i=1}^{n+1}}
(-1)^{n(i-1)}\tau_{(1\dots n+1)} ^{i-1}
\{\!\{D_i, \{\!\{D_{i+1}, \dots ,D_{n+1}, \dots , D_{i-1}\}\!\}_{ \mu_1 \dots  \mu_n}\}\!\}_L
\\
+ {\displaystyle \sum_{i=1}^{n­+1}} (-1)^{n(i+1)}
\tau_{(12\dots n+1)}^i \{\!\{ D_{i+1}, D_{i+2}, \dots ,D_{n+1}, D_1,  \dots ,\{\!\{D_{i-1}, D_i \}\!\}\}\!\}
_{ \mu_1 \dots  \mu_n , L} 
\end{array}$$
\end{theorem}

\begin{proof}

$$\Phi \left ( \underline{d}_\mathbb{L} (\mu_1 \dots \mu_n)\right )
(D_1 \otimes \dots \otimes D_{n+1})=
({\color{black}\underline{I}})+(\underline{II})$$
where 
$$\begin{array}{l}
({\color{black}\underline{I}})=\mu_n(D_{n+1})^{\prime \prime} D_1 \left (\mu_1 (D_2)^\prime \right ) \otimes 
\mu_1 (D_2)^{\prime \prime}\mu_2(D_3)^\prime \otimes \otimes \dots \otimes 
\mu_{n-1}(D_n)^{\prime \prime}\mu_n(D_{n+1})^\prime\\ 
+\dots + 
{\displaystyle \sum_{i=2}^n}(-1)^{i-1}
\mu_n(D_{n+1})^{\prime \prime}\mu_1(D_1)^\prime \otimes \dots \otimes 
D_i \left ( \mu_{i-1}(D_{i-1})^{\prime \prime} \mu_i (D_{i+1})^\prime\right ) \otimes \dots \otimes 
\mu_{n-1}(D_n)^{\prime \prime}\mu_n (D_{n+1})^\prime \\
(-1)^n  D_{n+1}\left ( \mu_n(D_n)^{\prime \prime}\right )^{\prime \prime}\mu_1(D_1)^\prime 
\otimes \mu_1(D_1)^{\prime \prime}\mu_2(D_2)^\prime \otimes \dots \otimes 
\mu_{n-1}(D_{n-1})^{\prime \prime}\mu_n(D_n)^\prime \otimes 
D_{n+1}\left ( \mu_n (D_n)^{\prime \prime}\right )^\prime\\
\end{array}$$ 
and 
$$\begin{array}{l}
(\underline{II})= 
{\displaystyle \sum_{i=1}^n}(-1)^i
\mu_n (D_{n+1})^{\prime \prime} \mu_1  (D_1)^\prime \otimes 
\mu_1(D_1)^{\prime \prime}\mu_2 (D_2)^\prime \otimes \dots \otimes \\
\mu_{i-1}(D_{i-1})^{\prime \prime}\mu_i \left (\{\!\{D_i, D_{i+1}\}\!\}_l^\prime \right )^\prime \otimes 
\{\!\{D_i, D_{i+1}\}\!\}_l^{\prime \prime} \otimes 
\mu_i\left (\{\!\{ D_i, D_{i+1}\}\!\}_l^\prime \right)^{\prime \prime} \mu_{i+1}(D_{i+2})^\prime \otimes \dots \otimes 
\mu_{n-1}(D_n)^{\prime \prime}\mu_n (D_{n+1})^\prime
\end{array}$$

For simplicity, we write $\sigma := \tau_{(1\dots n)}$ and $\tau := \tau_{(1\dots n+1)}$. 
Then 
$$\{\!\{-, \dots, -\}\!\}_{\mu_1 \dots \mu_n}={\displaystyle \sum_{i=1}^n}(-1)^{(n-1)(i-1)}
\mu_{\sigma^i(1)}\otimes \dots \otimes \mu_{\sigma^{i}(n)}.$$
Using the formulas of ${\color{black}\underline{I}}$ and $\underline{II}$, one writes  
$d_\mathbb{L} (\{\!\{-, \dots ,-\}\!\}_{\mu_1 \dots \mu_n})(D_1 \otimes \dots \otimes D_{n+1})=(I)+(II)$, 
In the computation of (I), the terms of the form 
$D_1()\otimes ...$ give 
$\{\!\{D_1, \{\!\{D_2, \dots ,D_{n+1}\}\!\}_{ \mu_1 \dots  \mu_n}\}\!\}_L$.  
More precisely :
$$\begin{array}{lll}
(I)& =& 
D_1 \left [ \mu_n (D_{n+1})^{\prime \prime} \mu_1 (D_2)^\prime\right ] \otimes 
\mu_1 (D_2)^{\prime \prime}\mu_2(D_3)^\prime \otimes \dots \otimes 
\mu_{n-1}(D_n)^{\prime \prime}\mu_n(D_{n+1})^\prime+ \dots\\
& + &
{\displaystyle \sum_{i=1}^{n} \sum_{j=1}^{n+1}}(-1)^{(n-1)(i-1)}(-1)^{n(j-1)}\tau^{j-1}
\left ( \mu_k \longleftarrow   \mu_{\sigma^{i-1}(k)}, D_l \longleftarrow D_{\tau^{j-1}(l)}\right )
\\
&=&
\{\!\{D_1, \{\!\{D_2, \dots ,D_{n+1}\}\!\}_{ \mu_1 \dots  \mu_n}\}\!\}_L + 
\dots  \\
&+&(-1)^{n(i-1)}\tau_{(1\dots n+1)} ^{i-1}
\{\!\{D_i, \{\!\{D_{i+1}, \dots ,D_{n+1}, \dots , D_{i-1}\}\!\}_{ \mu_1 \dots  \mu_n}\}\!\}_L
+ \dots  \\
&= &{\displaystyle \sum_{i=1}^{n+1}}
(-1)^{n(i-1)}\tau_{(1\dots n+1)} ^{i-1}
\{\!\{D_i, \{\!\{D_{i+1}, \dots ,D_{n+1}, \dots , D_{i-1}\}\!\}_{ \mu_1 \dots  \mu_n}\}\!\}_L
\end{array}$$
where the notation $\left ( \mu_k \longleftarrow   \mu_{\sigma^i(k)}, D_l \longleftarrow D_{\tau^j(l)}\right )$ means 
that we reproduce  the previous term,   replacing $\mu_k$ by $\mu_{\sigma^i(k)}$ and 
$D_l$ by $D_{\tau^j(l)}$. \\

When computing (II), the terms finishing by $\{\!\{D_n, D_{n+1}\}\!\}^{\prime \prime}$ give 
$$ (-1)^n\{\!\{D_1, \dots , D_{n-1}, \{\!\{D_n, D_{n+1}\}\!\}^\prime \}\!\}_{ \mu_1 \dots  \mu_n}
\otimes \{\!\{D_n, D_{n+1}\}\!\}^{\prime \prime}$$
More precisely :  
$$\begin{array}{l}
(II)= \\
(-1)^n \mu_n \left ( \{\!\{D_n, D_{n+1}\}\!\}^\prime\right )^{\prime \prime}\mu_1(D_1)^\prime  \otimes 
\mu_1 (D_1)^{\prime \prime} \mu_2 (D_2)^\prime \otimes \dots \otimes 
\mu_{n-1}(D_{n-1})^{\prime \prime}\mu_n \left ( \{\!\{D_n, D_{n+1}\}\!\}^\prime \right )\otimes 
\{\!\{ D_n, D_{n+1}\}\!\}^{\prime \prime}\\
+\dots + {\displaystyle \sum_{i=1}^{n} \sum_{j=1}^{n+1}}(-1)^{(n-1)(i-1)}(-1)^{n(j-1)}\tau^{j-1}
\left ( \mu_k \longleftarrow   \mu_{\sigma^{i-1}(k)}, D_l \longleftarrow D_{\tau^{j-1}(l)}\right )\\
= (-1)^n\{\!\{D_1, \dots , D_{n-1}, \{\!\{D_n, D_{n+1}\}\!\} \}\!\}_{ \mu_1 \dots  \mu_n}
+ \dots +(-1)^{2n}\tau^n  \left ( D_k \longleftarrow D_{\tau^{n+1}(D_k)} \right ) \\
= {\displaystyle \sum_{i=1}^{n+1}} (-1)^{n(i+1)}
\tau_{(12\dots n+1)}^i \{\!\{ D_{i+1}, D_{i+2}, \dots , \{\!\{D_{i-1}, D_i \}\!\}\}\!\}
_{ \mu_1 \dots  \mu_n , L} 
\end{array}$$
\end{proof}
 
 \begin{remark} 
 In the case where $A=k$, $\mathbb{L}=B$ is a finite dimensional  algebra and $\{\!\{a,b\}\!\}=ab-ba$, it is easy to see  (using theorem \ref{explicit formula}) that $\dfrac{T_A(\mathbb{L}^*)}{[T(\mathbb{L}^*), T(\mathbb{L}^*)]}$ is isomorphic to the cochains of the cyclic cohomology
 $$C^*_n=\{f: B^{\otimes n} \to k, \; f\circ \tau^{-1}_{(1\dots n)} =(-1)^{(n-1)}f\}$$
 and $d_{\mathbb L}$ is the differential of the Hochschild complex. Thus
  the complex 
  $\left ( \dfrac{T_A(\mathbb{L}^*)}{[T(\mathbb{L}^*), T(\mathbb{L}^*)]}, d_{\mathbb{L}}\right ) $ 
  computes the cyclic cohomology of $B$. \\\\
 
 \end{remark}

 We will now study more in detail the case where $A$ is a smooth double Poisson algebra with double bracket defined by the double biderivation  
 $P=\delta \Delta$ and $\mathbb{L}=\Omega^1_A$. Note that 
 $\Omega^1_A=\mathbb{D}er_k (A)^*$ and $\mathbb{D}er_k (A)=(\Omega_{A}^1)_*$ so that if $D$ is in 
 $\mathbb{D}er _k(A)$, then $\sigma D \in (\Omega^1_kA)^*$.\\

 \begin{proposition}
 
1) If $a,f \in A$, then $\underline{d}_\mathbb{L}(a)(df)=-\{\!\{f, \{P,a\} \}\!\}$ .
 
2)  If $D_1, \dots , D_n$ are in $\mathbb{L}$, then 
 $$\underline{d}_\mathbb{L}(\sigma D_1, \dots , \sigma D_n)(da_1 \otimes \dots \otimes da_{n+1}) = 
  - \{\!\{a_1, \{\!\{a_2, \dots \{\!\{a_{n+1}, \{P, D_1 \dots D_n \}\}\!\}_\mathbb{L} \dots \}\!\}_\mathbb{L}$$
  \end{proposition}
  
 \begin{proof} 
 Let us prove that $\underline{d}_\mathbb{L}(a)(df)=\{\!\{f, \{\delta \Delta ,a\} \}\!\}$.
  With our definition, we get 
 $$\begin{array}{rcl}
 \underline{d}_\mathbb{L}(a)(df)&=& 
 \{\!\{f,a\}\!\}_{\color{black}P}^\prime \otimes \{\!\{f,a \}\!\}_{\color{black} P}^{\prime \prime}\\
&=& \Delta(a) ^\prime \delta (f)^{\prime \prime} \otimes \delta (f)^\prime \Delta (a)^{\prime \prime}
-\delta (a)^\prime \Delta (f)^{\prime \prime}\otimes \Delta (f)^\prime \delta (a)^{\prime \prime}
\end{array}$$

On the other hand :
$$\begin{array}{rcl}
\{\delta \Delta , a\}&=& -\{\!\{a, \delta \Delta \}\!\}^{\prime \prime} \{\!\{a, \delta \Delta \}\!\}^\prime\\
&=& -\{\!\{a,\Delta\}\!\}^{\prime \prime}\delta \{\!\{a,  \Delta\}\!\}^\prime 
+ \{\!\{a, \delta\}\!\}^{\prime \prime} \Delta \{\!\{a, \delta \}\!\}^\prime\\
&=& \Delta (a)^\prime \delta \Delta (a)^{\prime \prime}-\delta (a)^\prime \Delta \delta (a)^{\prime \prime}
\end{array}$$
Hence 
$\{\delta \Delta ,a \}(df)= \delta (f)^\prime \Delta(a)^{\prime \prime}\otimes \Delta (a)^\prime \delta (f)^{\prime \prime}- \Delta (f)^\prime \delta(a)^{\prime \prime}\otimes \delta (a)^\prime \Delta (f)^{\prime \prime}$.\\
  
  The proof of 2) is a consequence of the two following lemmas. 
  
  \begin{lemma}
  $\underline{d}_\mathbb{L} (\sigma D) (da \otimes db)=-\{\!\{a, \{\!\{b, \{P,D\} \}\!\}_L$
  \end{lemma}
  
 \begin{proof}
  $$\begin{array}{l}
  \underline{d}_\mathbb{L}(\sigma D) (da \otimes db)\\
  =-D(a)^{\prime \prime}\otimes \{\!\{b, D(a)^\prime \}\!\}_P + 
 \{\!\{a, D(b)^{\prime \prime}\}\!\}_P \otimes D(b)^{\prime } -
  \tau_{(23)}\{\!\{\sigma D, \{\!\{a, \{P,b\} \}\!\} \}\!\}_L\\ 
 = -D(a)^{\prime \prime}\otimes \{\!\{b, D(a)^\prime \}\!\}_P + 
  \{\!\{a, D(b)^{\prime \prime}\}\!\}_P \otimes D(b)^{\prime } 
 - D \left ( \{\!\{a,b\}\!\}^\prime\right )^{\prime \prime} \otimes \{\!\{a,b\}\!\}^{\prime\prime } 
 \otimes D \left ( \{\!\{a,b\}\!\}^\prime\right )^{\prime }\\ 
 = -D(a)^{\prime \prime}\otimes \{\!\{b, D(a)^\prime \}\!\}_P + 
 \{\!\{a, D(b)^{\prime \prime}\}\!\}_P \otimes D(b)^{\prime } -
  \tau_{(132)}\{\!\{ D, \{\!\{a, \{P,b\} \}\!\} \}\!\}_L\\ 
= -D(a)^{\prime \prime}\otimes \{\!\{b, D(a)^\prime \}\!\}_P + 
 \{\!\{a, D(b)^{\prime \prime}\}\!\}_P \otimes D(b)^{\prime } -
  \{\!\{ a, \{\!\{\{P,b\}, D \}\!\} \}\!\}_L+   \tau_{(123)}\{\!\{ \{P,b\}, \{\!\{D, a \}\!\} \}\!\}_L\\ 
  =-D(a)^{\prime \prime}\otimes \{\!\{b, D(a)^\prime \}\!\}_P + 
 \{\!\{a, D(b)^{\prime \prime}\}\!\}_P \otimes D(b)^{\prime } -
  \{\!\{ a, \{\!\{\{P,b\}, D \}\!\} \}\!\}_L +   
 \tau_{(123)}
 {\color{black} \left [\{\!\{b,  \{\!\{D, a \}\!\}^\prime \}\!\}_P \otimes \{\!\{D, a \}\!\}^{\prime\prime} \right ]}\\ 
={\color{black}+}\{\!\{a, D(b)^{\prime \prime}\}\!\}_P \otimes D(b)^{\prime } +  \{\!\{ a, \{\!\{\{P,b\}, D \}\!\} \}\!\}_L \\
={\color{black}+} \{\!\{a, D(b)^{\prime \prime}\}\!\}_P \otimes D(b)^{\prime } + \{\!\{ a,  \{ P, \{\!\{b, D \}\!\}\} \}\!\}_L - 
 \{\!\{ a,  \{\!\{b, \{P, D \} \}\!\} \}\!\}_L \\
= -  \{\!\{ a,  \{\!\{b, \{P, D \}\}\!\} \}\!\}_L 
  \end{array}$$
  where we used the {\color{black}formula
$  \{\!\{a,b\}\!\}=\{\!\{a, \{P,b\}\}\!\}.$}
  \end{proof}
  
  \begin{lemma} If $\delta_1, \dots , \delta_n$ are in $\mathbb{D}er (A)$ and $a_1, \dots ,a_n \in A$, then 
  $$\begin{array}{l}
  \{\!\{a_1, \{\!\{a_2 , \dots ,\{\!\{a_n , \delta_1\dots \delta_n \}\!\}_L \}\!\}_L \dots \}\!\}_L
  = \\
  \delta_1 (a_1)^{\prime \prime} \otimes \delta_1(a_1)^\prime \delta_2(a_2)^{\prime \prime } \otimes
   \delta_2(a_2)^\prime \delta_3 (a_3)^{\prime \prime} \otimes \dots 
  \otimes \delta_{n-1}(a_{n-1})^\prime \delta_n (a_n)^{\prime \prime} \otimes \delta_n(a_n)^\prime
  \end{array}$$
  
  \end{lemma}
  
  Proof by induction on $n$. 
  
  \end{proof}
 
 \begin{proposition} Let $A$ be a smooth double Poisson algebra. 
  If $\mathbb{L}=\Omega^1_A$, the differential $d_\mathbb{L}$ coincides with the double Poisson cohomology defined by Pichereau and Van Weyer (\cite{Pichereau-vanWeyer}).
 \end{proposition}
 
 This proposition follows from the previous theorem. \\

 {\bf The Lie derivative :}
 
 The definition of  the Lie derivative is more complicated for a general double Lie algebroid than in the 
  case where  $\mathbb{L}=\mathbb{D}er(A)$ (\cite{CB-PE-VG}). \\
 
Let $D$ be an element of $\mathbb{L}$. The Lie derivative along $D$ is the map 
$L_D : T_A(\mathbb{L}^*) \to 
T_A(\mathbb{L}^*) \otimes  T_A(\mathbb{L}^*)$ defined by the following conditions : 
\begin{itemize}
\item $\underline{L}_D(a)=D(a)$
\item  If $\lambda \in \mathbb{L}^*$, 
$$\begin{array}{rcl}
\underline{L}_{D}(\lambda)(\Delta )&=&
\lambda (\Delta)^\prime \otimes D \left (\lambda (\Delta )^{\prime \prime} \right )
-\tau_{(23)}\left [ \lambda \left ( \{\!\{D, \Delta \}\!\}_r^{\prime \prime}\right ) \otimes  \{\!\{D, \Delta \}\!\}_r^{\prime }\right ]\\
&+&
D\left ( \lambda (\Delta)^\prime\right ) \otimes \lambda (\Delta )^{\prime \prime} 
-\tau_{(32)} \left [\lambda \left ( \{\!\{D, \Delta \}\!\}_l^{ \prime }\right ) \otimes 
 \{\!\{D, \Delta \}\!\}_l^{\prime \prime} \right ]\\
\end{array}$$
$\underline{L}_D (\lambda)\in A \otimes \mathbb{L}^* \oplus \mathbb{L}^* \otimes A$ as the map 
$$\Delta  \mapsto \lambda (\Delta )^\prime \otimes D\left ( \lambda (\Delta)^{\prime \prime} \right ) -
\tau_{(32)} \left [\lambda \left (\{\!\{D, \Delta \}\!\}_r^{\prime\prime}\right ) \otimes 
\{\!\{D, \Delta \}\}_r^{\prime } \right ]$$
belongs to $\mathbb{L}^*\otimes A$ and the map 
$$\Delta  \mapsto D\left (\lambda (\Delta)^\prime \right )\otimes \lambda (\Delta)^{\prime \prime}
- \tau_{(32)} \left [
\lambda ( \{\!\{D, \Delta \}\!\}^\prime_l) \otimes \{\!\{D, \Delta \}\!\}_l^{\prime \prime})\right ]$$
is  in $A\otimes \mathbb{L}^*$

\item $\underline{L}_D$ is a degree preserving double derivation from $T_A(\mathbb{L}^*)$ to 
$T_A(\mathbb{L}^*) \otimes  T_A(\mathbb{L}^*)$.

\end{itemize}



 Standard formulas involving the differential 
 $\underline{d}_{\mathbb{L}}$, 
 contraction and Lie derivative hold in the case where $\mathbb{L}=\mathbb{D}er(A )$ 
 (\cite{CB-PE-VG}, \cite{VdB1}). We prove now that they hold for any 
  double Lie algebroids. In particular, one has the Cartan identity :
 
 \begin{proposition} \label{unreduced} The map $\underline{d}_\mathbb{L}$ is a degree one derivation of 
 $T_A(\mathbb{L}^*)$ and can be extended to a degree one derivation of 
 $T_A(\mathbb{L}^*)\otimes T_A(\mathbb{L}^*)$ as follows :
 $$\forall (\alpha , \beta) \in T_A(\mathbb{L}^*), \quad 
 \underline{d}_\mathbb{L}(\alpha \otimes \beta )=\underline{d}_{\mathbb{L}}(\alpha)\otimes \beta +
 (-1)^{\vert \alpha \vert } \alpha \otimes \underline{d}_\mathbb{L}(\beta).$$
 
 One has the following properties : For any $D$ and $\Delta$ in $\mathbb{L}$ :

1) $ \underline{d}_\mathbb{L} \circ \underline{i}_{D}+\underline{i}_{D}\circ \underline{d}_\mathbb{L } 
=\underline{L}_D$

2) $\{\!\{\underline{i}_D, \underline{i}_\Delta \}\!\}_l=
 \{\!\{\underline{i}_D, \underline{i}_\Delta \}\!\}_r=0$.\\
 
 \end{proposition}
 
\begin{remark} Consequences of these formulas will be seen later. \end{remark}

\begin{proof}

2) It is enough to prove the relation 
$\underline{d}_\mathbb{L} \circ \underline{i}_{D}+\underline{i}_{D}\circ \underline{d}_\mathbb{L} 
=\underline{L}_D$
on elements of $A$ and $\mathbb{L}^*$. 

On elements of $A$, it is obvious. On elements $\lambda \in \mathbb{L}^*$, 
we give the main steps of the computation ~:
$$\begin{array}{rcl}
(\underline{i}_D\circ \underline{d}_\mathbb{L} + \underline{d}_\mathbb{L} \circ \underline{i}_D)(\lambda )(\Delta )
&=& 
\underline{d}_\mathbb{L}\lambda (D, \Delta ) -\underline{d}_\mathbb{L}\lambda (\Delta ,D)+
\underline{d}_\mathbb{L}\left ( \lambda (D)\right )(\Delta)\\
&=&\underline{d}_\mathbb{L}\lambda (D, \Delta ) -\underline{d}_\mathbb{L}\lambda (\Delta ,D)+ 
\underline{d}_\mathbb{L}\left ( \lambda (D)^\prime \right )\otimes \lambda (D)^{\prime \prime}+
\lambda (D)^\prime \otimes \underline{d}_\mathbb{L}\left ( \lambda (D)^{\prime \prime} \right )\\
&=&\underline{L}_D(\lambda )(\Delta )
\end{array}$$

3) The relations $\{\!\{\underline{i}_D, \underline{i}_\Delta \}\!\}_l^\prime=0$ and 
$ \{\!\{\underline{i}_D, \underline{i}_\Delta \}\!\}_r^{\prime \prime}=0$ are easy to check on elements of $A$ and 
$\mathbb{L}$.
\end{proof}


\section{From double to classical} 

We now explain how to go from the double picture to the classical picture by the use of representation spaces (\cite{VdB1}). 
In the section, we assume by simplicity that $B=k$. \\

Let $A$ be a  $k$-algebra an $N\in \N$. Denote by $Rep_N(A)$ the representation space 
$Hom_{alg} \left (A, M_N(k) \right )$
The coordinate ring of $Rep_N(A)$ is 
$$ {\mathcal O}_N(A):= 
k[Rep_N(A)]=\dfrac{k\left [a_{p,q}, (p,q) \in [1,N], a \in A \right ]}{<a_{p,q}b_{q,r} -(ab)_{p,r}>}$$
For any element $x \in Rep_N(A)$ one has  
$$a_{p,q}(x)=x(a)_{p,q}$$
\begin{examples}
1) If $A=\dfrac{k[t]}{(t^n)}$, then $Rep_N(A)=\{Q\in M_{N}(k) \mid Q^n=(0)\}.$

2) If $A=k<x_1, \dots , x_q>$, then $Rep_N(A)=M_N(k)\oplus \dots \oplus M_N(k).$
\end{examples}

\begin{theorem}(\cite{VdB1}) 
If $(A, \{\!\{-,-\}\!\}$)
 is a double Poisson algebra, then ${\mathcal O}_N(A)$, endowed with  the bracket determined by
 $$\{a_{i,j},b_{u,v}\}=\{\!\{a,b\}\!\}^\prime_{u,j}\{\!\{a,b\}\!\}^{\prime \prime}_{i,v}$$
 is a Poisson algebra. 
 \end{theorem}

If $a\in A$, one introduces $X(a)$ the $M_n(k)$ valued function on $Rep_N(A)$ defined by 
$X(a)_{i,j}=a_{i,j}$. One has the relation $X(ab)=X(a)X(b)$ and one  defines $Tr: \dfrac{A}{[A,A]} \to 
{\mathcal O}_N(A) , \quad \overline{a}\to Tr(X(a))={\displaystyle \sum_i }a_{i,i }$. 

$GL_N(k)$ acts by conjugation on $Rep_N(A)$. Using a result of Lebruyn and Procesi
\cite{L-P},  the following theorem is shown in (\cite{CB1})

\begin{theorem} 
If the $k$-algebra $A$ is finitely generated, the map 
$$ \begin{array}{rcl}
Tr: \dfrac{A}{[A,A]} &\to &
{\mathcal O}_N(A)^{GL_N(k)} \\
 \overline{a} &\to &Tr(X(a))={\displaystyle \sum_i }a_{i,i }
\end{array}$$
is an isomorphism of Lie algebras. 
\end{theorem}

If  $M$ is an  $A^e$-module, one defines (\cite{VdB2}) the ${\mathcal O}_N(A)$-module $(M)_N$
$$ (M)_N=\dfrac{k\left [m_{i,j} , \; m \in M\right ]}
{<a_{i,u}m_{u,j}-(am)_{i,j}, a_{u,j}m_{i,u}-(ma)_{i,j}>}$$
If $m\in M$, one introduces $Tr(m)={\displaystyle \sum_{i=1}^N m_{i,i} \in (M)_N}$.

We will be interested in the case where $M=\mathbb{L}$ is a double Lie Rinehart algebra. \\

If $\delta \in \mathbb{D}er (A)$, one  defines (\cite {VdB1}) the corresponding derivation  on 
${\mathcal O}_N(A)$ by 
$$\delta_{i,j}(a_{u,v})=\delta (a)^\prime_{u,j}\delta(a)^{\prime \prime}_{i,v}.$$ If $\delta=\delta_1 \dots \delta_n$, one sets 
$\delta_{i,j}=\delta_{1,ia_1}\delta_{2,a_1a_2}\dots \delta_{n, a_{n-1}j} \in \wedge 
Der ({\mathcal O}_N(A))$. In other words, this can be rewritten by the relation 
$X(\delta)=X(\delta_1)\dots X(\delta_n)$. 

\begin{proposition}(\cite{VdB1})  If $P,Q \in DA$, then the following relation holds : 
$$\{P_{i,j}, Q_{u,v}\}=\{\!\{P,Q\}\!\}^\prime_{u,j}\{\!\{P,Q\}\!\}_{i,v}^{\prime \prime}$$
 where $\{\!\{-,-\}\!\}$ denotes the Schouten bracket on $DA$ and $\{-,-\}$ the Schouten bracket between poly-vector fields on $Rep_N(A)$. 
\end{proposition}

\begin{proposition} The Trace map $Tr$ (\cite{VdB1})
$$\begin{array}{rcl}
Tr \;:\dfrac{DA}{[DA,DA]}& \to & \wedge Der \left ( {\mathcal O}_{N}(A)\right )\\
\delta & \mapsto & Tr X(\delta)
\end{array}$$
is a Lie algebra homomorphism if both side are equipped with the Schouten bracket. 
\end{proposition}

\begin{theorem}(\cite{VdB2})

There exists a unique Lie  algebroid structure  
$\left ( (\mathbb{L})_N, {\mathcal O}_N(A) , [ , ], \omega \right )$ with bracket  
$[, ]$ and anchor $\omega$ determined below by the equalities below :
$$\begin{array}{l}
\omega (X_{i,j})(a_{u,v})=\omega (X)(a)_{u,j}^\prime \omega (X)(a)_{i,v}^{\prime \prime}\\
\left [X_{i,j}, Y_{u,v} \right ] = \{\!\{X,Y \}\!\}^\prime _{u,j}  \{\!\{X,Y \}\!\}^{\prime \prime} _{i,v}
\end{array}$$
\end{theorem}

\begin{proof} The only thing that is not obvious is the Jacobi identity and the fact that $\omega$ is a Lie algebra morphism. We need to prove the identities 
$$\begin{array}{l}
\left [ X_{i,j}, [Y_{u,v},Z_{k,m}]\right ]+\left [ Y_{u,v}, [Z_{k,m},X_{i,j}]\right ] +
\left [ Z_{k,m}, [X_{i,j},Y_{u,v}]\right ]=0\\
\omega \left ([X_{i,j}, Y_{u,v} \right ])= \left [ \omega(X_{i,j}), \omega (Y_{u,v})\right ]
\end{array}$$
These two identities follows from the double Jacobi identities by a straightforward computation.\end{proof}


If $D\in\mathbb{L}$, we consider   the matrix  $X(D)=\left (D_{i,j} \right )$ as being  with values in ${\color{black}\bigwedge}_{{\mathcal O}_N(A)}(\mathbb{L})_N$ and we set 
$$\forall D_1 \otimes \dots \otimes D_n, \quad X(D_1 \dots D_n)=
X(D_1)\dots X(D_n) \in {\color{black}\bigwedge}_{{\mathcal O}_N(A)} (\mathbb{L})_N.$$
 One defines the trace map (as in \cite{VdB1}) by 
 $Tr : T_A(\mathbb{L}) \to {\color{black}\bigwedge}_{\mathcal O} (\mathbb{L})_N $  by 
 $Tr(D_1\otimes \dots \otimes D_n)=TrX(D)$.
\begin{proposition} 
 If $D$ and $\Delta$ are in $\mathbb{L}$, one has the following equality 
$[Tr (D), Tr (\Delta ) ]=Tr \left ( \{\!\{D, \Delta  \} \!\}\right )$ where the left hand side  involves the Schouten bracket on ${\color{black}\bigwedge}_{{\mathcal O}_N(A)} (\mathbb{L})_N$ 
and the right hand side  involves the double Schouten bracket.
\end{proposition}
\begin{proof}  it is a straightforward computation. \end{proof}

Let $M$ be an $A^e$-module. To the $A^e$-module  
$$M^*=\{ \lambda : M \to A \otimes A \mid 
\lambda (amb)=\lambda (m)^\prime b \otimes  a \lambda (m)^{\prime \prime}\},$$
 one associates 
$$(M^*)_N=\dfrac{k\left [\lambda_{i,j}, \; \lambda \in M^* \right ]}
{<a_{i,u}\lambda_{u,j}-(a\lambda)_{i,j}, a_{u,j}\lambda_{i,u}-(\lambda a)_{i,j}>}$$
If $\lambda \in M^*$ and $m \in M$, $\lambda_{i,j}$ defines an element of 
$Hom_{\mathcal{O}_N(A)}\left ((M)_N,\mathcal{O}_N(A) \right )$ by 
$$\lambda_{i,j}(m_{u,v})=\lambda (m)^\prime_{i,v}\lambda (m)^{\prime \prime}_{u,j}$$

For us,  $M$ will be a double Lie--Rinehart algebra $\mathbb{L}$ or its dual ${\mathbb{L}^*}$. 

\begin{remark} If $M=\mathbb{D}er(A)$, one recovers the equality $(\underline{d}a)_{i,j}=da_{i,j}$  (\cite{VdB1}).\end{remark}

We go on mimicking Van den Berg's construction :\\

If $\Lambda=\lambda_1 \otimes \lambda_2 \otimes \dots \otimes \lambda_n \in T^n (\mathbb {M}^*)$, 
 we define 
$$\Lambda_{i,j}= \lambda_{1,ii_1}\lambda_{2,i_1i_2}\dots \lambda_{n,i_{n-1}j}
 \in {\color{black}\bigwedge} \left (  (\mathbb{M})_N^*\right )$$
and $X(\Lambda)=(\Lambda_{i,j})_{i,j}$. The latter   is a matrix with values in 
${\color{black}\bigwedge} \left (  \mathbb{M})_N^*\right )$. 

\begin{lemma}
Let us identify $\Lambda = \lambda_1 \otimes \lambda_2 \otimes \dots \otimes \lambda_n$ to 
$\Phi \in Hom_{A^e}( \mathbb{L}^{\otimes n}, A^{\otimes n+1})$. If we write  
  $$\Phi (D_1\otimes \dots \otimes D_n)=
\Phi^{(1)}(D_1\otimes \dots \otimes D_n)\otimes \dots \otimes \Phi^{(n+1)}(D_1\otimes \dots \otimes D_n),$$
then 
$$\begin{array}{rcl}
\Phi _{i,j}\left (D_{1,u_1v_1} \otimes \dots \otimes D_{n, u_nv_n} \right )&=&
\displaystyle \sum_{\sigma \in S_n} (-1)^n
\Phi^{(1)}(D_{\sigma (1)}\otimes \dots \otimes D_{\sigma (n)})_{iv_1} 
\Phi^{(2)}(D_{\sigma (1)} \otimes \dots \otimes D_{\sigma (n)})_{u_1v_2} \dots \\
&&\Phi^{(n-1)}(D_{\sigma (1)}\otimes \dots \otimes D_{\sigma (n)})_{u_{n-1}v_n}
\Phi^{(n)}(D_{\sigma (1)}\otimes  \dots \otimes D_{\sigma (n)})_{u_{n}j}
\end{array}$$
\end{lemma}

\begin{proof}
$$\begin{array}{l}
{\color{black}\Phi^{(1)}}(D_{\sigma (1)} \dots D_{\sigma (n)})_{iv_1} 
\Phi^{(2)}(D_{\sigma (1)} \dots D_{\sigma (n)})_{u_1v_2} 
\Phi^{(n-1)}(D_{\sigma (1)} \dots D_{\sigma (n)})_{u_{n-1}v_n}
\Phi^{(n)}(D_{\sigma (1)} \dots D_{\sigma (n)})_{u_{n}j}\\
=
\lambda_1 \left ( D_{\sigma (1)}\right )^\prime_{iv_1}
\left [\lambda_1 (D_{\sigma (1)})^{\prime \prime} \lambda_2(D_{\sigma (2)})^\prime \right ]_{u_1 v_2} \dots
\left [\lambda_{n-1}(D_{\sigma (n-1)}) ^{\prime \prime} \lambda_n (D_{\sigma (n)})^\prime \right ]_{u_n ,v_n}
\lambda(D_{\sigma (n)})^{\prime \prime}_{u_n,j}\\
= \lambda_{1,ia_1}(D_{\sigma (1), u_1v_1}) \lambda_{2,a_1a_2}(D_{\sigma (2), u_2v_2})\dots
 \lambda_{n,a_{n-1}j}(D_{\sigma (n), u_nv_n})\\
= (\lambda_1\otimes \lambda_2 \otimes  \dots \otimes \lambda_n)_{i,j}
 \left (D_{\sigma(1)  u_1v_1}\otimes \dots \otimes D_{\sigma (n), u_nv_n} \right )
 \end{array}$$
 \end{proof}

\begin{lemma}
(i) If $a\in A$, then $d_\mathbb{L}(a)_{i,j}=d_{(\mathbb{L})_N}(a_{i,j})$.

(ii) If $\lambda \in \mathbb{L}^*$, then $\underline{d}_\mathbb{L}(\lambda )_{i,j}=
d_{(\mathbb{L})_N}(\lambda_{i,j})$
\end{lemma}

\begin{proof}

Let us prove (i). On one hand, one has :

$$
\underline{d}_\mathbb{L}(a)_{i,j}(D_{u,v})=
\underline{d}_\mathbb{L}(a)(D)^\prime_{i,v} \underline{d}_\mathbb{L}(a)(D)_{u,j}^{\prime \prime}=
D(a)^\prime_{i,v}D(a)^{\prime \prime}_{u,j}$$
On the other hand :

$$\begin{array}{rcl}
d_{(\mathbb{L})_N}(a_{i,j})(D_{u,v})=\omega (D_{u,v})(a_{i,j})=D(a)^\prime_{i,v}D(a)^{\prime \prime}_{u,j}.
\end{array}$$

Let us now prove (ii). In the following computation $\left ( D\leftrightarrow \Delta , u\leftrightarrow k , v\leftrightarrow p\right )$ means the same expression as before exchanging $D$ with $\Delta$, $u$ with $k$, $v$ with $p$. 
$$\begin{array}{rcl}
\underline{d}_\mathbb{L} (\lambda )_{i,j}(D_{k,p},\Delta_{u,v})
&=& D\left ( \lambda (\Delta) ^\prime\right )^\prime _{i,p}
D\left ( \lambda (\Delta) ^\prime\right )^{\prime \prime} _{k,v}
 \lambda (\Delta )^{\prime \prime} _{u,j}-
\lambda (D)^\prime_{i,p}  \Delta \left ( \lambda (D)^{\prime \prime}\right )^\prime _{k,v}
 \Delta \left ( \lambda (D)^{\prime \prime}\right )^{\prime\prime}_{u,j}\\
&-&\lambda \left (\{\!\{D,\Delta \}\!\}^\prime_l\right )^\prime_{i,p}  (\{\!\{D, \Delta \}\!\}^{\prime\prime}_{l})_{k,v}
\lambda \left (\{\!\{D,\Delta \}\!\}^\prime_l\right )^{\prime \prime}_{u,j}\\
&-&\left ( D\leftrightarrow \Delta , u\leftrightarrow k , v\leftrightarrow p\right )\\
&=& D_{k,p}\left ( \lambda (\Delta) ^\prime_{i,v}\right )
 \lambda (\Delta )^{\prime \prime} _{u,j}-
\lambda (D)^\prime_{i,p}  \Delta_{u,v} \left ( \lambda (D)^{\prime \prime}_{k,j}\right )^\prime \\
&-&\lambda_{i,j} \left ((\{\!\{D,\Delta \}\!\}^\prime_{l})_{ u,p}\right )  (\{\!\{D, \Delta \}\!\}^{\prime\prime}_{l})_{,k,v}\\
&-&\left ( D\leftrightarrow \Delta , u\leftrightarrow k , v\leftrightarrow p\right )\\
&=& D_{k,p}\left ( \lambda_{i,j}(D_{u,v})\right ) -
\Delta_{u,v}\left (\lambda_{i,j}(D_{k,p}) \right )\\
&-& \lambda_{i,j} \left ((\{\!\{D,\Delta \}\!\}^\prime_{l})_{u,p}\right ) ( \{\!\{D, \Delta \}\!\}^{\prime\prime}_{l})_{k,v}
- \lambda_{i,j} \left ((\{\!\{D,\Delta \}\!\}^{\prime\prime}_{r})_{ k,v}\right )  
(\{\!\{D, \Delta \}\!\}^{\prime}_{r})_{u,p}\\
&=& D_{k,p}\left ( \lambda_{i,j}(D_{u,v})\right ) -
\Delta_{u,v}\left (\lambda_{i,j}(D_{k,p}) \right ) 
- \lambda_{i,j} \left ([D_{k,p},\Delta_{u,v}]\right )\\
&=& d_{(\mathbb{L)}_N}(\lambda_{i,j})(D_{k,p}, \Delta_{u,v})

\end{array}$$ \end{proof}


\begin{proposition}
 Let $Tr$ be the Trace map : if $\Phi \in T_A(\mathbb{L}^*)$, 
 $Tr (\Phi)=Tr \left [ X(\Phi) \right ]$

a) If $\lambda  \in \mathbb{L}^*$ and $D \in \mathbb{L}$, 
then $Tr \left ( \lambda (D)\right )= Tr (\lambda )\left ( Tr(D) \right )$.

b) If $\Phi \in T_A(\mathbb{L}^*)$, one has $Tr d_\mathbb{L}(\Phi)=d_{(\mathbb{L})_N}\left ( (Tr \Phi) \right )$
\end{proposition}

\begin{proof} a) is an easy computation :

$$Tr (\lambda ) \left ( Tr D) \right)= {\displaystyle \sum_{i,j}}\lambda_{i,i}(D_{j,j})
=\lambda (D)^\prime_{i,j} \lambda (D)^{\prime \prime}_{j,i}=Tr \left ( \lambda (D) \right ).$$

Let us now prove b). 

Let us first prove by induction on the degree of $\Phi$ that 
$\underline{d}_{\mathbb{L}} (\Phi)_{i,j}=d_{(\mathbb{L})_N}(\Phi_{i,j})$. 

For $deg(\Phi)=1$, we have already proved it in a previous lemma.

Assume that it is proved for $\deg (\Phi)=n$ and let us prove it 
for $\Phi \lambda$ if $\lambda \in \mathbb{L}^*$.

$$\begin{array}{rcl}
\underline{d}_\mathbb{L} (\Phi \lambda)_{i,j}&=& 
\left ( \underline{d}_\mathbb{L}(\Phi)\lambda + \Phi \underline{d}(\lambda)\right )_{i,j}\\
&=& \underline{d}_\mathbb{L}(\Phi)_{i,i_1}\lambda_{i_1,j} + 
\Phi _{i,i_1}\underline{d}_\mathbb{L}(\lambda)_{i_1 ,j}\\
&=&d_{(\mathbb{L})_N} (\Phi_{i,i_1}) \lambda_{i_1,j}+ \Phi_{i,i_1}d_{(\mathbb{L})_N}(\lambda_{i_1,j})\\
&=&d_{(\mathbb{L})_N}\left ( (\Phi \lambda )_{i,j}\right )
\end{array}$$
Then $Tr\left [ d_\mathbb{L}(\Phi)\right ]=
{\displaystyle \sum_i} d_\mathbb{L}(\Phi )_{i,i} =
{\displaystyle \sum_i} d_{(\mathbb{L})_N}(\Phi_{i,i})=d_{(\mathbb{L})_N}(Tr(\Phi))$. 
\end{proof}

\begin{theorem}
The Tr induces a map from 
$\dfrac{T_A(\mathbb{L}^*)}{[T_A(\mathbb{L}^*), T_A (\mathbb{L}^*)]}$ to 
$\bigwedge_{{\mathcal O}_N(A)} (\mathbb{L}^*)_N$ that sends $d_\mathbb{L}$ to 
$d_{(\mathbb{L})_N}$ 
\end{theorem}


\section{Reduced contraction and Lie derivative}

In this section, we generalize results of \cite{CB-PE-VG}. 

\begin{definition}
Let $\mathbb{L}$ be a double Lie algebroid. If $\Theta$ is in $\mathbb{L}$, 
one defines the reduced Lie derivative and the reduced contraction by :
$$\begin{array}{l}
\iota_\Theta : T_A^n({\mathbb{L}^*}) \to T_A^{n-1}(\mathbb{L}^*), 
\quad \alpha \mapsto {}^\circ (\underline{i}_{\Theta}\alpha )\\
{\mathcal L}_\Theta : T_A^n(\mathbb{L}^*) \to T_A^{n}(\mathbb{L}^*), 
\quad \alpha \mapsto  {}^\circ(\underline{L}_{\Theta}\alpha  )\\
\end{array}$$
\end{definition}

Explicitely, if $\alpha_1, \alpha_2, \dots , \alpha_n$ are in $\mathbb{L}^*$, one has 
$$\begin{array}{l}
\iota_\Theta (\alpha_1 \alpha_2 \dots \alpha_n)=
{\displaystyle \sum_{k=1}^n (-1)^{k(n-k+1)}}
\alpha_k (\Theta)^{\prime \prime}\cdot \alpha_{k+1}\dots \alpha_n \alpha_1 \dots \alpha_{k-1} \cdot
\alpha_k (\Theta)^\prime \\
{\mathcal{L} }_\Theta (\alpha_1 \alpha_2 \dots \alpha_n)=
{\displaystyle \sum_{k=1}^n (-1)^{k(n-k+1)}}
{\color{black}\underline{L}}_{\Theta}(\alpha_k )^{\prime \prime}\cdot \alpha_{k+1}\dots \alpha_n \alpha_1 \dots \alpha_{k-1} 
\cdot {\color{black} \underline{L}}_{\Theta}(\alpha_k )^\prime \\
\end{array}$$


\begin{proposition}
1) For any $\Theta \in \mathbb{L}$,  the following equalities of endomorphisms of $T_A(\mathbb{L}^*)$ hold :
$$\underline{d}_\mathbb{L} \circ \iota_\Theta +\iota_\Theta \circ \underline{d}_\mathbb{L} =
{\mathcal L}_\Theta ,\quad 
\underline{d}_\mathbb{L}\circ {\mathcal L}_\Theta ={\mathcal L}_\Theta \circ \underline{d}_\mathbb{L} $$

2) The maps $\underline{d}_\mathbb{L}$, ${\mathcal L}_\Theta$ and $\iota_\Theta$ descend to maps from 
$DR^\bullet (\mathbb{L})$ to $DR^\bullet(\mathbb{L})$ denoted respectively 
$d_\mathbb{L}$, ${\mathcal L}_\Theta$ and $\iota_\Theta$.  One has 
$$d_\mathbb{L} \circ \iota_\Theta +\iota_\Theta \circ d_\mathbb{L} =
{\mathcal L}_\Theta ,\quad 
d_\mathbb{L}\circ {\mathcal L}_\Theta ={\mathcal L}_\Theta \circ d_\mathbb{L} $$

3) For any $\delta$ and $\Delta$ in $\mathbb{L}$, one has 
$\underline{i}_D \iota_{\Delta}+\sigma_{(12)}\underline{i}_\Delta \iota_\delta =0$ 
(as maps from $T_A(\mathbb{L}^*)$ to 
$T_A(\mathbb{L}^*) \otimes T_A(\mathbb{L}^*)$).
\end{proposition}

\begin{remark} The previous proposition is proved for 
$\mathbb{L}=\mathbb{D}er (A)$ in \cite{CB-PE-VG} but our proof is different. \end{remark}
\begin{proof}  

1)  follows by applying $^\circ()$ to  the relation 
$$\underline{L}_\Theta =
\underline{d}_\mathbb{L} \circ \underline{i}_\Theta  + \underline{i}_\Theta \circ \underline{d}_\mathbb{L}$$
and using proposition \ref{reduced and unreduced}. 

 


2) followed from proposition \ref{reduced and unreduced} and from proposition \ref{unreduced}.
  
\end{proof}






\end{document}